\documentclass[14pt]{amsart}
\usepackage{cases}
\usepackage{amsmath}
\usepackage{amsfonts}
\usepackage{bm}
\usepackage{color}
\usepackage{amssymb}
\usepackage{xypic}
\usepackage[all]{xy}
\usepackage{mathrsfs}
\usepackage{amsthm}
\usepackage[all]{xy}
\usepackage{graphics}
\usepackage{array}
\usepackage{graphicx} 
\usepackage{epstopdf}
\usepackage{float}
\usepackage{tikz}
\usepackage{tikz-cd}
\usetikzlibrary{arrows}
\usetikzlibrary{graphs}
\usepackage{tikz}
\usepackage{mathtools}
\usepackage{hyperref}
\usepackage{delarray}

\usepackage{subcaption}
\usepackage{booktabs}
\usepackage{extarrows}
\usepackage{rotating}

\newtheorem{thm}{Theorem}[section]
\newtheorem{lem}[thm]{Lemma}
\newtheorem{defi}[thm]{Definition}
\newtheorem{cor}[thm]{Corollary}
\newtheorem{prop}[thm]{Proposition}
\newtheorem{ex}[thm]{Example}
\newtheorem{rmk}[thm]{Remark}

\title
{Finest positroid subdivisions from maximal weakly separated collections}

\author{Gleb A. Koshevoy \and Fang Li \and Lujun Zhang*}

\address{Gleb A. Koshevoy
\newline
Institute for Information Transmission Problems,
Russian Academy of Science of Moscow,
 Moscow 127051, Russia.}
\email{koshevoyga@gmail.com}

\address{Fang Li
\newline School of Mathematical Sciences,
Zhejiang University,
Yuhangtang Road 866,
Hangzhou, Zhejiang 310058,
China P.R.}
\email{fangli@zju.edu.cn}

\address{Lujun Zhang
	\newline
	Thr corresponding author. School of Mathematical Sciences,
	Zhejiang University,
	Yuhangtang Road 866,
	Hangzhou, Zhejiang 310058,
	China P.R. }
\email{12135007@zju.edu.cn}

\newcommand{\lra}{\longrightarrow}

\newcommand{\ra}{\rightarrow}
\newcommand{\sdp}{\times\kern-.2em\vrule height1.1ex depth-.05ex}
\newcommand{\epi}{\lra \kern-.8em\ra}

\makeatletter
\@addtoreset{equation}{section}
\makeatother

 \normalbaselines
\begin{document}
\renewcommand{\thefootnote}{\alph{footnote}}
\setcounter{footnote}{-1} \footnote{Keywords: positive tropical Grassmannian, positroid subdivision, weak separation, plabic graph, cluster algebra.

The corresponding author: Lujun Zhang, Email:12135007@zju.edu.cn.}

\renewcommand{\thefootnote}{\alph{footnote}}
\maketitle
\bigskip
	
\begin{abstract}
We adopt a formal and algebraic approach of Early \cite{E2} to study the positive tropical Grassmannian $\operatorname{Trop}^+ Gr_{k,n}$.
 Specifically, we deal with positroid subdivision of hypersimplex induced by translated blades from any maximal weakly separated collection. One of our main results gives a necessary and sufficient condition on a maximal weakly separated collection to form a positroid subdivision of a hypersimplex corresponding to a simplicial cone
in $\rm Trop^+Gr_{k,n}$. For k = 2 our condition says that any weakly separated
collection of two-elements sets gives such a simplicial cone, and all cones are
of such a form.

 We also show that the maximality of any weakly separated collection is preserved
under the boundary map, which armatively answers a question by Early in \cite{E1}. Plabic graphs, invented by Postnikov \cite{P}, are of use
in proving this result. As a corollary, we get 
that all those positroid subdivisions are the finest. Thus, the flip of two maximal weakly separated
collections corresponds to a pair of adjacent maximal cones in positive tropical Grassmannian.

\noindent\textbf{2020 Mathemathics Subject Classification: 05E45.}

\end{abstract}


\tableofcontents

\vspace{0.5cm}


\vspace{1cm}
\section{Introduction and notations}
For a positive integer $n$, let $[n]$ denote the ordered set of elements $\{1, 2, \cdots, n\}$.

For two sets $A$ and $B\subset [n]$ we say (i) $A$ is  {\em below} $B$ if,  for any $i\in A\setminus B$ and $j\in B\setminus A$, we have $i<j$;
(ii) $A$ {\em  splits} $B$ if both $A\setminus B$ and $B\setminus A$ are nonempty, and   $B\setminus A$ can be expressed as a disjoint union $B'\sqcup B''$  of nonempty
subsets so that  $B'$ is below $A\setminus B$ and  $A\setminus B$ is below $B''$.

Sets $A$, $B\subset [n]$  are called {\em weakly separated} (from each other) if either $A$ is below $B$, or $B$ is below $A$, or $A$ splits $B$ and $|A|\ge |B|$, or $B$ splits
$A$ and $|B|\ge |A|$.  A collection $ W\subset 2^{[n]}$ is called {\em weakly separated } ($w$-collections for short) if any two of its members are weakly separated.

These notions were introduced by Leclerc and Zelevinsky in \cite{LZ} where their importance is demonstrated, in particular, in connection with the problem of characterizing
quasicommuting quantum flag minors of a quantum matrix. 
They conjectured that all inclusion-wise maximal collections of this sort have the same
cardinality. Danilov et al in \cite{DKK}  answered affirmatively
this conjecture for flag varieties and Grassmanians.

Cluster algebras, invented  by Sergey Fomin and Andrei Zelevinsky in \cite{FZ}  are commutative algebras whose generators and relations are constructed in a recursive manner. Among these algebras, there are the algebras of
homogeneous coordinates on the Grassmannians, on the flag varieties and on many
other varieties which play an important role in geometry and representation theory.

From construction of flips between $w$-collections in  \cite{DKK} follows that  such collections are in bijection with Pl\"ucker cluster seeds of the cluster algebra of  $\mathbb C[GL_n/N]$, where $N$ is unipotent radical of $GL_n(\mathbb C)$, whose cluster variables are minors. Maximal $w$-collections of subsets $[n]$ of size $k$ are in bijection with Pl\"ucker cluster seeds of the cluster algebra of  $\mathbb C[Gr_{k,n}]$, the coordinate ring of Grassmaniian of $k$-planes, whose cluster variables are minors.

Early in \cite{E1} established  a connection between weakly separated collections and positroid subdivisions of hypersimplex $\Delta_{k,n}$ of the unit cube, which is defined as the section of the unit cube given by $\Delta_{k,n}=[0,1]^n\cap \{\sum\limits_{i=1}^n x_i=k\}$, $1<k<n$.

Regular positroid subdivisions of  $\Delta_{k,n}$ are important for understanding combinatorics of the positive tropical Grassmanian  $\rm Trop^+Gr_{k,n}$.  Namely, the maximal cones of $\rm Trop^+Gr_{k,n}$ correspond to the finest positroid subdivisions, whereas its rays correspond to the coarsest.


The positive tropical Grassmannian $\rm Trop^+Gr_{k,n}$ was introduced by Speyer and Williams in  \cite{SW}, defined as the space of realizable positive tropical linear spaces. 
In recent years, considerable  research has been carried out to explore its connections with other areas of mathematics and physics
(see for example \cite{SW2}, \cite{LPW}, \cite{E2}).

Regular positroid subdivisions arise as affinity  areas of tropical Pl\"ucker functions (TP functions for short, which are also called positive tropical vectors) defined on $\Delta_{k,n}$. The later form a subsclass of $M^{\#}$-functions invented by Murota (see \cite{M}). On one hand side,  for any maximal $w$-collection $W$ of $k$-sets, any TP-function  on $\Delta_{k,n}$ is defined by their values at the  vertices of $\Delta_{k,n}$ corresponding to $W$. However, such a relation  between TP-functions and weakly separated collections is not useful here, because TP-functions are not stable under summation and even under convolution ($M^{\#}$-functions are stable under convolutions).

On the other hand, Early in \cite{E1} used the vertices of $\Delta_{k,n}$ which are labelled by  a weakly separated collection to translate the standard blade. A standard blade $\beta$ is defined as the union of codimension $1$ faces of a complete fan in the hyperplane $H_0:=\{x \in \mathbb{R}^n\;|\;\sum_{i=1}^{n}x_i=0\}$. Translating $\beta$ to a vertex $e_J$ of $\Delta_{k,n}$ yields a translated blade $\beta_J$, which induces a multi-split of the hypersimplex.
 One of main results in \cite{E1} is that, for a weakly separated collection $W$,   the collection of translated blades $\{\beta_J,\, J\in W\}$ gives a positroid subdivision of $\Delta_{k,n}$. 
   Recall that a $k$-split (for $k \geq 2$) is a polyhedral subdivision into $k$ maximal cells intersecting along a common inner face of codimension $k-1$. The following theorem due to Herrmann guarantees that the positroid subdivision induced by single blade $\beta_J$ corresponds to a ray in $\rm Trop^+Gr_{k,n}$.

\vspace{0.3cm}
 \begin{thm}(\cite{H})
     A $k$-split is a coarsest regular subdivision.
 \end{thm}
  In \cite{E2} the notion of blade arrangement is extended to weighted blade arrangement which is a $\mathbb{R}$-linear combination of translated blades. This generalization allows to implement the positive tropical Grassmannian by weighted blade arrangements that satisfy compatibility condition and positivity condition under the boundary map $\partial_L$ (Definition \ref{def4}). Here we use a slightly different setting in defining this collection of weighted blade arrangements $\mathcal{Z}_{k,n}$, and prove that $\mathcal{Z}_{k,n}$ is exactly $\rm Trop^+Gr_{k,n}$ (see Lemma \ref{lem4}). Let $\overline{\mathcal{Z}_{k,n}}$ be the quotient fan of $\mathcal{Z}_{k,n}$ modulo the $n$-dimensional linearity space $\mathcal{F}_{k,n}:=\rm span \{\beta_I\;|\;I\;is\;frozen\}$. When $k=2$, for a nonfrozen maximal $w$-collection $W$, the polyhedron cone $\mathcal P_W:=\{\sum_{J\in W}\mathbb R_{\ge 0}\beta_J\}$ is simplicial of maximal dimension  in $\overline{\mathcal Z_{2,n}}$. However, when $k>2$, this is not always true.

For $k>2$, we define the set $\mathscr{S}$ of nonfrozen $w$-collections using the boundary map $\partial_L$ (for details see Section \ref{subsec1}). The collections of  $\mathscr{S}$ are  partially ordered by inclusion. The following theorem provides a necessary and sufficient condition for  a collection $W$, such that the cone $\mathcal P_{W}$ is  simplicial of maximal dimension in $\overline{\mathcal Z_{k,n}}$.
\vspace{0.3cm}
  \begin{thm}(Theorem \ref{thm1})
  \label{thm5}
For any nonfrozen maximal $w$-collection $W$, $\mathcal{P}_{W}$ is maximal simplicial in $\overline{\mathcal{Z}_{k,n}}$ if and only if $W$ is minimal in $\mathscr{S}$.
  \end{thm}
  
  For better understanding elements of $\mathscr{S}$, stability of maximality of the weak separation under the boundary map is an important question,  a problem posed in \cite{E1}.  We have the following result.
  \vspace{0.2cm}
  \begin{thm}(Theorem \ref{thm2})
      Let $W$ be a maximal $w$-collection in $\binom{[n]}{k}$, then $\partial_j(W)$ is also a maximal $w$-collection in $\binom{[n]\backslash \{j\}}{k-1}$.
  \end{thm}
Our proof of this theorem uses  Postnikov's theory of plabic graphs, a combinatorial framework for studying the totally nonnegative Grassmannian \cite{P}. From these main results we get the following corollary. 


\vspace{0.2cm}
\begin{cor}(Corollary \ref{cor5})
    The translated blades of a maximal weakly separated collection induce a finest regular positroid subdivision of $\Delta_{k,n}$. 
\end{cor}


\vspace{0.3cm}
This paper is organized as follows.

 We begin with introducing a new construction of positive tropical Grassmannian following \cite{E1, E2} in Section \ref{sec2}. 
Then in Section \ref{subsec1} we prove our first theorem. In Section \ref{sec3}, we introduce a boundary map on a reduced plabic graph and prove several useful lemmas. Then we give the proof of the second theorem and an example in Section \ref{subsec2}.

We use the following notations throughout the paper.

 $\binom{[n]}{k}$ denote the set of $k$-subsets of $[n]$. 
 We write $[a,b]$ for the closed cyclic interval from $a$ to $b$ and $(a,b):=[a,b]\backslash \{a,b\}$ for the open cyclic interval. If $L$ is a subset of $[n]$ and $a \in [n]\setminus L$, we abbreviate $L \cup \{a\}$ by $La$. For any $I \in \binom{[n]}{k}$, we use $e_I:=\sum_{i \in I} e_i$ and $x_I:\sum_{i \in I}x_i$ to denote the sum of unit vectors and coordinates. If no confusion, we write the set $I=\{i_1,i_2,\cdots,i_k\}$ by $i_1i_2\cdots i_k$ . 
 \vspace{1cm}
\section{Positive tropical Grassmannian and the first theorem}
\label{sec2}
\subsection{Preliminaries on positive tropical Grassmannian and positroid subdivisions}\quad


 \begin{defi}
For positive integers $k <n$, the \textbf{non-negative real Grassmannian} $Gr_{k,n}^{\geq 0}(\mathbb{R})$ is defined to be the quotient $Gr_{k,n}^{\geq 0}(\mathbb{R})=GL_k^{>0}(\mathbb{R})\backslash Mat_{k,n}^{\geq 0}(\mathbb{R})$, where $GL_k^{>0}(\mathbb{R})$ is the group of $k \times k$-matrices with positive determinants and  $Mat_{k,n}^{\geq 0}(\mathbb{R})$ denotes the set of all row-full-rank $k \times n$ real matrices $V$ such that all $k \times k$ minors $\{p_I(V)\;|\;I \in \binom{[n]}{k}$ are nonnegative.
\end{defi}
\vspace{0.2cm}

\begin{defi}
    For $\mathcal{M} \subset \binom{[n]}{k}$, let
 \[
 S_{\mathcal{M}}:=\{V \in Gr_{k,n}^{\geq 0}(\mathbb{R})\;|\;p_I(V) >0 \;if \;and \;only \; if \; I \in \mathcal{M}\}.
 \]
 If $S_{\mathcal{M}} \neq \emptyset$, then $\mathcal{M}$ is referred to as a \textbf{positroid}, and $S_{\mathcal{M}}$ as a \textbf{positroid cell}. The associated \textbf{positroid polytope} is defined as $P_{\mathcal{M}} := \operatorname{Conv} \{ e_J \mid J \in \mathcal{M}\}$, where $\operatorname{Conv}\{\cdot\}$ denotes the convex hull of the points. A polyhedral subdivision is called a \textbf{positroid subdivision} if every cell in the subdivision is a positroid polytope.
\end{defi}

\vspace{0.2cm}

In \cite{SW}, Speyer and Williams defined the positive tropical Grassmannian using tropicalizing the positive part of the Grassmannian over the ring of Puiseux series. Subsequently in \cite{SW2}, they proved that the positive tropical Grassmannian coincides with the positive Dressian, which is the collection of all positive tropical Pl\"ucker vectors.
\begin{defi}
The  \textbf{positive tropical Grassmannian} $\rm Trop^{+}Gr_{k,n}$ is the set of points $p=(p_I)_{I \in  \binom{[n]}{k}} \in \mathbb{R}^{\binom{n}{k}}$ such that
\begin{equation}
\label{E1}
p_{Lac}+p_{Lbd}={\rm min}\left\{p_{Lab}+p_{Lcd},\;p_{Lbc}+p_{Lad}\right\}
\end{equation}
for any $1 \leq a <b<c<d \leq n$ and $L \in \binom{[n]\backslash \{a,b,c,d\}}{k-2}$ and $p$ is called a \textbf{positive tropical Pl\"{u}cker vector}.
\end{defi}

\vspace{0.3cm}

\vspace{0.2cm}
Therefore, the positive tropical Grassmannian can be related to the regular subdivision of a certain polytope called hypersimplex. A \textbf{hypersimplex}  $\Delta_{k,n}$ 
is the convex hull of points $\left\{e_I\;|\;I \in \binom{[n]}{k}\right\}$.

If we take $p = (p_I) \in \mathbb{R}^{\binom{n}{k}}$, place the value $p_I$ at the vertex $e_I$ of $\Delta_{k,n}$, and consider the convex hull $\operatorname{Conv}\{(e_I,p_I)\;|\;I \in \binom{[n]}{k}\}$ which forms a polyhedron in $\mathbb{R}^{n+1}$. Projecting the lower faces of this polyhedron back onto $\Delta_{k,n}$ yields a polyhedral subdivision. Such a vector $p$ is called a \textbf{height function} and the resulting subdivision is called a \textbf{regular subdivision}, denoted by $\mathcal{D}_p$


\vspace{0.3cm}
\begin{thm}(\cite{LPW}, Theorem 9.12)
$p \in \rm Trop^{+}Gr_{k,n}$ if and only if $\mathcal{D}_p$ is a positroid subdivision of $\Delta_{k,n}$.
\end{thm}
Thus $\rm Trop^{+}Gr_{k,n}$ is equipped with a \textbf{secondary fan structure} such that $p$ and $p'$ are in the same relative interior of a cone if and only if they induce the same positroid subdivision of  $\Delta_{k,n}$.
Elements in a ray of $\rm Trop^{+}Gr_{k,n}$, when used as height functions, induce the \textbf{coarsest positroid subdivision}, whereas elements in the relative interior of a maximal cone induce the \textbf{finest positroid subdivision}. On the other hand, relations (\ref{E1}), being a set of linear equalities and inequalities, naturally defines the positive tropical Grassmannian $\rm Trop^{+}Gr_{k,n}$ as a polyhedral fan. This fan is called the \textbf{Pl\"ucker fan}.

\vspace{0.3cm}
\begin{thm}(\cite{OPS2}, Theorem 4.1)
    The Pl\"ucker fan structure coincides with the secondary fan structure in $\rm Trop^{+}Gr_{k,n}$.
\end{thm}

\vspace{0.2cm}

In the following sections of this paper, we will prove that Early's maximal blade arrangements generate finest positroid subdivisions of $\Delta_{k,n}$. The following theorem by Speyer and Williams (2021) provides equivalent characterizations of the finest positroid subdivisions.

\begin{thm}(\cite{SW2}, Theorem 6.6)
\label{thm3}
   Let $p$ be a positive tropical Pl\"ucker vector. Then the following statements are equivalent,
    \begin{enumerate}
    \item $\mathcal D_p$ induces a finest subdivision of $\Delta_{k,n}$, that is, for any positive Pl\"ucker vector $p'$ such that $p+p'$ is positive Pl\"ucker,  there holds $\mathcal D_p=\mathcal D_{p+p'}$.
    
\item Every octohedral face ${Lab,Lbc,Lcd,Lad,Lac,Lbd}$ in  $\mathcal D_p$ is subdivided in two pyramids, each of which contains the non-separated diagonal ${Lac,Lbd}$.

\item $p_{Lab}+p_{Lcd} \neq p_{Lad}+p_{Lbc}$ for any $a<b<c<d$ and $L \in \binom{[n]\backslash \{a,b,c,d\}}{k-2}$.

    \end{enumerate}
\end{thm}

\vspace{1cm}
\subsection{Blade arrangement and weighted blade arrangement}\quad

In this and the next section, we will review Early's blade arrangement model, as well as the open question regarding boundary maps proposed by Early in \cite{E1}.

\vspace{0.5cm}
\begin{defi}

A \textbf{decorated ordered set partition} of type $\Delta_{k,n}$ is an ordered set partition $(S_1,S_2,\cdots,S_l)$ together with a list of nonnegative integers $(s_1, \cdots, s_l)$, where
\begin{enumerate}
    \item $\sum\limits_{i=1}^ls_i=k$,
    \item $0 \leq s_i \leq |S_i|-1$ for each $i=1,2,\cdots,l$.
\end{enumerate}
We denote this decorated partition by $\{(S_1)_{s_1},(S_2)_{s_2},\cdots,(S_l)_{s_l}\}$ and the collection of decorated ordered set partitions of type $\Delta_{k,n}$ by $OSP(\Delta_{k,n})$.
\end{defi}

Given a $\{(S_1)_{s_1},(S_2)_{s_2},\cdots,(S_l)_{s_l}\} \in OSP(\Delta_{k,n})$, we use $[(S_1)_{s_1}, (S_2)_{s_2},\cdots, (S_l)_{s_l}]$ to denote the (translated) polyhedral cone in $H_{k,n}$  
 formulated by the following facet-defining inequalities
\begin{equation}
\label{E2}
\begin{split}
    x_{S_1} &\geq s_1,\\
    x_{S_1 \cup S_2} &\geq s_1+s_2,\\
    &\vdots\\
    x_{S_1 \cup \cdots \cup S_{l-1}} &\geq s_1+ \cdots s_{l-1}.
\end{split}
\end{equation}

\begin{prop}(\cite{E1})
    For a partition in $OSP(\Delta_{k,n})$, the $l$ cyclically shifted polyhedral cones $[(S_1)_{s_1}, (S_2)_{s_2},\cdots, (S_l)_{s_l}]$, $[(S_l)_{s_l}, (S_1)_{s_1},\cdots, (S_{l-1})_{s_{l-1}}]$, $\cdots$, $[(S_2)_{s_2}, (S_3)_{s_3},\cdots, (S_1)_{s_1}]$ form a complete simplicial fan in $H_{k,n}$.
\end{prop}

\vspace{0.3cm}
\begin{defi}
    The \textbf{blade} corresponding to the partition, denoted by $((S_1)_{s_1},(S_2)_{s_2},\cdots,(S_l)_{s_l})$, is the union of the codimension-1 faces of its associated cyclic fan in the affine space $H_{k,n}$. Formally, if $\partial P$ denotes the boundary of a closed cone P, then
    \[
    ((S_1)_{s_1},(S_2)_{s_2},\cdots,(S_l)_{s_l})=\bigcup\limits_{j=1}^{l}\partial[(S_j)_{s_j},(S_{j+1})_{s_{j+1}},\cdots,(S_{j-1})_{s_{j-1}}].
    \]
    In particular, when $S_i=\{i\}$ and $s_i=0$ for all the $i \in [n]$, the blade is called \textbf{standard blade} denoted by $\beta:=((1,2,\cdots,n))$.
\end{defi}
The  standard blade $\beta = ((1,2,\cdots,n))$ is characterized by the following union,
\[
\beta=((1,2,\cdots,n))=\bigcup\limits_{j=1}^{n}\partial[j,j+1,\cdots,j-1]
\]
where $\Pi_j :=[j,j+1,\cdots,j-1]=\left\{t_1(e_1-e_2)+t_2(e_2-e_3)+ \cdots t_j\widehat{(e_j-e_{j+1})}+\cdots t_n(e_n-e_1)\right\}$ ($j+1:=1$ if $j=n$). Let the \textbf{translated blade} $\beta_J=((1,2,\cdots,n))_{e_{J}}$ be the translation of $\beta$ from origin to the vertex $e_J$ of $\Delta_{k,n}$. 

The blade can be realized as the locus where a piecewise linear function attains its minimum value at least twice, i.e., as a \textbf{tropical hypersurface}. Particularly, for standard blade $\beta$, this piecewise-linear function $h(x)$ on hyperplane $H_{0,n}=\{x \in \mathbb{R}^n\;|\; \sum\limits_{i=0}^n x_i=0\}$ can be given by
\[
h(x)={\rm min} \{L_1(x),L_2(x),\cdots,L_n(x)\}
\]
where 
\[
L_i(x)=x_{i+1}+2x_{i+2}+\cdots+(n-i)x_n+x_1+\cdots+(n-1)x_{i-1}
\]

\vspace{0.3cm}
\begin{prop}(\cite{E1})
\label{prop5}
    The affine areas of $h(x)$ are exactly $\Pi_1,\Pi_2,\cdots,\Pi_n$, and 
    \[
    {\rm min} \{L_1(x),L_2(x),\cdots,L_n(x)\}=L_i(x)
    \]
    for $x \in \Pi_i$. Thus, the standard blade $\beta$ is a tropical hypersurface.
\end{prop}

\vspace{0.3cm}
Now since $\Pi_1,\Pi_2,\cdots,\Pi_n$ form a complete simplicial fan in $H_{0,n}$, then for any $x \in (H_{0,n} \cap \mathbb{Z}^n)$, there exists a unique maximal subset $\{i_1,i_2,\cdots,i_l\} \subset [n]$ such that $x \in \cap_{j \in \{i_1,i_2,\cdots,i_l\}}\Pi_j$. Moreover, it can be verified that $\{e_j-e_{j+1}\;|\;j \in \{i_1,i_2,\cdots,i_l\}^c\}$ is a minimal system of generators for the monoid $\left(\cap_{j \in \{i_1,i_2,\cdots,i_l\}}\Pi_j\right) \cap \mathbb{Z}^n$. Thus we obtain a unique \textbf{positive expression}
\[
x=\sum\limits_{j \in \{i_1,i_2,\cdots,i_l\}^c}t_j(e_j-e_{j+1})
\]
where $t_j >0$ for $j \notin \{i_1,i_2,\cdots,i_l\}$. The set $\{i_1,i_2,\cdots,i_l\}^c$ is called the \textbf{support} of $x$, denoted by $Supp(x)$. In particular, for $x = e_J - e_I$ with any distinct $k$-subsets $I$ and $J$, the positive integer $\sum_{j \in {i_1,i_2,\cdots,i_l}^c} t_j$ is called the \textbf{distance}  from $e_I$ to $e_J$, denoted by $d(e_I,e_J)$. Clearly, $d(e_I,e_J)+d(e_J,e_I)=n$ whenever $I$ and $J$ satisfy $|I \setminus J| = |J \setminus I| = 1$.
\vspace{0.3cm}

\begin{lem}
    \label{lem3}
    Let $h(x-e_J)$ be the translation of $h(x)$ from the origin to $e_J$, then 
    \[
    h(e_I-e_J)=-d(e_J,e_I)
    \]
\end{lem}
\begin{proof}
    We first take the unique positive expression $e_I-e_J=\sum_{j \in Supp(e_I-e_J)}t_j(e_j-e_{j+1})$. Without loss of generality, suppose $e_I-e_J \in \Pi_1$ (In other words, $t_n=0$), then
    \[
    \begin{split}
        h(e_I-e_J) &=L_1(e_I-e_J)\\
        &=(e_I-e_J) \cdot (0,1,2,\cdots,n-1)\\
        &=(\sum_{j \in Supp(e_I-e_J)}t_j(e_j-e_{j+1})) \cdot (0,1,2,\cdots,n-1)\\
        &=-\sum_{j \in Supp(e_I-e_J)}t_j
    \end{split}
    \]
\end{proof}

For a given $J \in \binom{[n]}{k}$, a decorated ordered set partition $\{(S_1){s_1}, (S_2){s_2}, \cdots, (S_l){s_l}\}$ can be constructed as follows. First, decompose $J$ into a disjoint union of cyclic intervals: $J = \bigsqcup\limits_{i=1}^{l} J_i$. Let $C_i$ be the cyclic interval connecting $J_{i-1}$ ($J_0:=J_l$) and $J_i$, such that $J_{i-1} \sqcup C_i \sqcup J_i$ itself forms a cyclic interval. Then, by defining $S_i = C_i \sqcup J_i$ and $s_i = |J_i|$ for $i = 1, 2, \cdots, l$ such that $1 \in S_1$, we obtain an element of $OPS(\Delta_{k,n})$.
\vspace{0.3cm}
\begin{thm}(\cite{E1})
    \label{thm6}
    Let $e_J$ be a vertex of $\Delta_{k,n}$, the translated blade $\beta_{J}=((1,2,\cdots,n))_{e_{J}}$ induces a $l$-split of $\Delta_{k,n}$ such that 
    \[
    ((1,2,\cdots,n))_{e_{J}} \cap \Delta_{k,n}=((S_1)_{s_{1}},(S_2)_{s_2},\cdots,(S_l)_{s_l}) \cap \Delta_{k,n},
    \]
    where $S_i=C_i \sqcup I_i $ and $s_i=|J_i|$ for $i=1,2, \cdots, l$ and $1 \in S_1$
\end{thm}
\vspace{0.2cm}

\begin{ex}
    Let $e_J$ be a vertex of $\Delta_{5,12}$ where $J=\{1,3,4,7,9\}$. Then $\beta_J$ induces a 4-split of $\Delta_{5,12}$ with 
   \begin{align*}
J_1 = \{1\},\quad & S_1 = \{10,11,12,1\}, & J_2 = \{3,4\},\quad & S_2 = \{2,3,4\}, \\
J_3 = \{7\},\quad & S_3 = \{5,6,7\}, & J_4 = \{9\},\quad & S_4 = \{8,9\}.
\end{align*}
\end{ex}

\vspace{0.3cm}
\begin{rmk}

    In Theorem~\ref{thm6}, the blade $\beta_J$ induces a decomposition of the hypersimplex $\Delta_{k,n}$ into $l$ full-dimensional polytopes
\[
    \widetilde{SM}_j := \left[ (S_j)_{s_j}, (S_{j+1})_{s_{j+1}}, \dots, (S_{j-1})_{s_{j-1}} \right] \bigcap \Delta_{k,n}, \quad \text{for } j = 1, 2, \dots, l.
\]
These polytopes meet along a unique common face of codimension~$l$, which is given by the product of simplices
\[
    \Delta_{s_1, S_1} \times \Delta_{s_2, S_2} \times \cdots \times \Delta_{s_l, S_l},
\]
where each $\Delta_{s_j, S_j} := \left\{ x \in [0,1]^{S_j} \;\middle|\; \sum_{i \in S_j} x_i = s_j \right\}$. Moreover, all the inequalities in \ref{E2} are facet-defining inequalities of $\widetilde{SM}_j$,

Following Oh~\cite{O}, each $\widetilde{SM}_j$ is called a \textbf{cyclic shifted dual Schubert matroid} and has been proved to be a positroid polytope. The subdivision induced by $\beta_J$ is trivial if and only if the index set $J$ is a cyclic interval. In this case, we refer to $\beta_J$ as \textbf{frozen}; otherwise, it is called \textbf{nonfrozen}.
\end{rmk}


\vspace{0.3cm}
\begin{defi}
    A \textbf{blade arrangement}  is a collection of translated blades $\{\beta_{J_1},\beta_{J_2}, \cdots, \beta_{J_l}\}$. A \textbf{weighted blade arrangement} is a formal $\mathbb{R}$-linear combination of some translated blades $\{\beta_{J_1},\beta_{J_2}, \cdots, \beta_{J_l}\}$.
\end{defi}

\vspace{0.3cm}
\begin{rmk}
    \label{rmk1}
    
    A blade $\beta_I$ can be regarded as a indicator function on $\Delta_{k,n}$, that is,
    \[
    \beta_I(x)=\begin{cases}
        1, &x \in \beta_I \cap \Delta_{k,n},\\
        0, &x \in \Delta_{k,n} \backslash \beta_I.
    \end{cases}
    \]
    A weighted blade arrangement $\sum\limits_Ic_I\beta_I$  can be regarded as a weighted sum of these indicator functions, i.e., for $x \in \Delta_{k,n}$,
    \[
    \sum\limits_{I}c_I\beta_I(x)=\sum\limits_{x \in \beta_I}c_I.
    \]
    In particular, we can take the blade arrangement as a weighted blade arrangement with positive coefficients.
\end{rmk}

For $k$-sets, weakly separation coincides with the chord separation.

\begin{defi}
    Two sets $I,J \in \binom{[n]}{k}$ are \textbf{weakly separated} if there do not exist $a,b,c,d \in [n]$, cyclically ordered, such  that
        $a,c \in I \backslash J$ and $b,d \in J\backslash I$.
    A collection of k-sets $W$ is called a \textbf{w-collection} if any two sets of $W$ are weakly separated. 
    $W$ is \textbf{maximally weakly separated} if for all $k$-sets $I \notin W$, $S \cup \{I\}$ is not weakly separated.
\end{defi}

\vspace{0.3cm}
In \cite{E1},  Early treated weak separation phenomena from the viewpoint of discrete geometry. We will see in Lemma \ref{lem4} that this is equivalent to telling whether two rays $\mathbb{R}_{\geq 0}\beta_I$ and $\mathbb{R}_{\geq 0}\beta_J$ are in some maximal cone of the positive tropical Grassmannian $\rm Trop^+Gr_{k,n}$.\\
\begin{thm}(\cite{E1})
    \label{thm4}
    The refinement of subdivisions $\{\beta_{J_1},\beta_{J_2}, \cdots, \beta_{J_l}\}$ is a positroid subdivision of $\Delta_{k,n}$ if and only if $\{J_1,J_2,\cdots,J_l\}$ is a weakly separated collection.
\end{thm}

\vspace{1cm}
\subsection{Construction of positive tropical Grassmannian by weighted blade arrangement}\quad\\


In this section, we explore the positive tropical Grassmannian using a more formal and algebraic approach, distinct from the geometric perspective employed by Speyer and Williams.
The advantage of this method is the ability to use a tool known as a boundary map, which allows for a recursive process descending to lower-dimensional $\Delta_{k,n}$ faces. Several open questions within this framework are raised in \cite{E1}. In the following sections, we will answer some of them.

For any subset $S \subseteq \mathbb{R}^n$, let $\partial_j(S)$ be the image of the set $S \cap \{ x_j = 1 \}$ under the canonical projection $(x_1, \dots, x_n) \mapsto (x_1, \dots, x_{j-1}, x_{j+1}, \dots, x_n)$. For example, $\partial_j(\Delta_{k,n}):=\Delta_{k-1,[n]\backslash \{j\}}=\{x \in [0,1]^{[n]\backslash \{j\}}\;|\;\sum\limits_{i \in {[n]\backslash \{j\}} }x_i=k-1\} $
\begin{lem}(\cite{E1})
   \label{lem5}
   
    Let $J=\{j_1,j_2,\cdots,j_k\} \in \binom{[n]}{k}$ with $1 \leq j_1<j_2<\cdots<j_k \leq n$, there holds
    \[
    \partial_j(\beta_J \cap \Delta_{k,n})=\beta_{J'}^{\{j\}} \cap \partial_j(\Delta_{k,n}),
    \]
    where $J=\{j_1,\cdots,j_k\}$ and $J'=J\backslash \{j_{a+1}\}$ if $j$ satisfies $j_{a} < j \leq j_{a+1}$ in the cyclic order. $\beta_{J'}^{\{j\}}:=((1,2,\cdots,\widehat{j},\cdots,n))_{e_{J'}}$ is the translated standard blade by deleting the $j$-th element
\end{lem}

Lemma \ref{lem5} motivates the following definition.

\begin{defi}
    We adopt the notations of Lemma \ref{lem5}. The \textbf{boundary map} $\partial_j$ on a set and translated blade is defined, respectively, as $\partial_j(J)=J'$ and $\partial_j(\beta_J)=\beta_{J'}^{\{j\}}$.
\end{defi}
We can extend the boundary map by replacing the singleton set $\{j\}$ with a subset $L=\{a_1, \cdots ,a_t\} \subset [n]$ satisfying $1 \leq t \leq k-2$.  The corresponding composite boundary map is defined as
    \[
    \partial_L:=\partial_{a_1} \circ \partial_{a_1} \circ \cdots \circ \partial_{a_t}
    \]
    which can be shown to be independent of the ordering of the elements $a_1, \cdots, a_t$.
    
Building on this, we consider the general setting. Let $\beta_J^{(L)}$ be the translated blade on $\partial_L(\Delta_{k,n})$ for a subset $L \subset [n]$ with $0 \leq |L| \leq k-2$ and $J \subset \binom{[n]\backslash L}{k-|L|}$.
 The boundary map $\partial_j$ can be defined  recursively on $\beta_J^{(L)}$ by,
\[
\partial_j(\beta_J^{(L)})=\left\{
\begin{aligned}
    &0  &j \in L\\
    &\beta_{J\backslash \{j_{a+1}\}}^{(L\cup \{j\})}  & j \notin L
\end{aligned}
\right.
\]
where $J=\{j_1, \cdots, J_{k-|L|}\}$ and $j_a < j \leq j_{a+1}$ in the  cyclic order of $[n]\backslash L$ inherited from $[n]$. The boundary map on a collection of blade arrangements can be extended linearly to weighted blade arrangement.

\begin{defi}
\label{def2}
    Let
    \[
    \mathfrak{B}_{k,n}^{\bullet}:= \bigoplus\limits^{n-(k-2)}_{m=0}\left(\bigoplus\limits_{L \in \binom{[n]}{m}}\mathfrak{B}_{k,n}^{(L)}\right),
    \]
    where $\mathfrak{B}_{k,n}^{(L)}$ is a vector space of formal linear span of $\beta_J^{(L)}$
    \[
    \mathfrak{B}_{k,n}^{(L)}:=\mathrm{span} \left\{\beta_J^{(L)}\;:\;J \in \binom{[n]\backslash L}{k-|L|}\right\}
    \]
    Let $\mathfrak{B}_{k,n}:=\mathfrak{B}_{k,n}^{(\emptyset)}$ be the top component of $\mathfrak{B}_{k,n}^{\bullet}$. Then $(\mathfrak{B}_{k,n}^{\bullet},\partial)$ is called the \textbf{hypersimplicial blade complex} where $\partial:=\partial_1+\partial_2+\cdots+\partial_n$ is the sum of all the boundary maps
    \[
    \partial:\bigoplus\limits_{L \in \binom{[n]}{m}}\mathfrak{B}_{k,n}^{(L)} \longrightarrow \bigoplus\limits_{L \in \binom{[n]}{m-1}}\mathfrak{B}_{k,n}^{(L)}
    \]
\end{defi}

\vspace{0.3cm}
\begin{rmk}
 The pair $(\mathfrak{B}_{k,n}^{\bullet},\partial)$ is not a chain complex in the usual sense, since $\partial^2 \neq 0$ even though $\partial_j^2 = 0$ for each $j \in [n]$.

 We will see further that a weighted blade corresponds to a function that admits a positroid decomposition of $\Delta_{k,n}$. However, the sum of weighted blades does not necessarily yield such a function, that is, it may no longer induce a positroid decomposition  of $\Delta_{k,n}$.

For this reason, we will later introduce a subset $\mathcal{Z}_{k,n}$ consisting of those formal sums of weighted blades that actually induce a positroidal decomposition.
\end{rmk}

\vspace{0.2cm}
 We consider a weighted blade arrangement $\mathcal{L}_J$ defined by
    \[
    \mathcal{L}_J:=\sum\limits_{e_{J_M} \in C_{J^{\bullet}}}(-1)^{1+|M|}\beta_{J_M},
    \]
    where we define the set
    \[
    C_{J^{\bullet}}=\left\{e_{J_M}=e_J+\sum\limits_{i=1}^{t}\delta_M(j_i)(e_{j_i-1}-e_{j_i})\;|\; M \subset J^{\bullet} \right\}.
    \]
    Here, $J^{\bullet}=\{j_1,\cdots,j_t\}$ is the set of initial elements of cyclic intervals of $J$ and $\delta_M$ is the characteristic function of subset $M \subset J^{\bullet}$ (i.e., $\delta_M(j_i) = 1$ if $j_i \in M$, and $\delta_M(j_i) = 0$ otherwise). Analogously, we use $J_{\bullet}$ to denote the set of final elements of cyclic intervals of $J$.
\vspace{0.3cm}

\begin{ex}
    For $k=7$, $n=15$ and $J=[1,2] \cup [4,5] \cup[9,11]$, we have $J^{\bullet}=\{1,4,9\}$ and 
    \[
    \mathcal{L}_J=-(\beta_{149}+\beta_{138}+\beta_{39,15}+\beta_{48,15})+(\beta_{38,15}+\beta_{49,15}+\beta_{139}+\beta_{148}).
    \]
    To prevent readers from misinterpreting notations like "15" as ${1} \cup {5}$, we separate it from the other numbers with a comma.
\end{ex}

\vspace{0.2cm}
\begin{prop}
\label{prop4}
We have the following relation
\[
\sum\limits_{I \in \binom{[n]}{k}} d(e_J,e_I)\mathcal{L}_I=n\beta_J-(\beta_{[1,k]}+\beta_{[2,k+1]}+\cdots+\beta_{[n,k-1]})
\]
where $d(e_J,e_I)$ is the smallest  positive number of steps from $e_J$ to $e_I$, where each step has to be one of the directions  $e_1-e_2,e_2-e_3,\cdots,e_n-e_1$. 
\end{prop}

\begin{proof}
    Let $K_{\bullet}=\{k_1',k_2',\cdots,k_t'\}$ be the set of final elements of cyclic intervals of $K$ for some $K \in \binom{[n]}{k}$ and $C_{K_{\bullet}}=\left\{e_{K_M}=e_K-\sum\limits_{i=1}^{t}\delta_M(k_i')(e_{k_i'}-e_{k_i'+1})\;|\; M \subset K_{\bullet} \right\}$. Then we calculate the coefficient of $\beta_K$ on the left side of the equation. It is equal to 
    \[
    \sum\limits_{M \subset K_{\bullet}} (-1)^{1+|M|}d(e_J,e_{K_M}).
    \]
    Let $e_K-e_J=\sum\limits_{j=1}^{n}t_j(e_j-e_{j+1})$ be the unique positive expression of $e_K-e_J$, so $d(e_J,e_K)=\sum\limits_{j=1}^{n}t_j$ and
    \[
    \begin{split}
        e_{K_M}-e_J &=e_{K_M}-e_K+e_K-e_J,\\
        &=-\sum\limits_{i=1}^{t}\delta_M(k_i')(e_{k_i'}-e_{k_i'+1})+\sum\limits_{j=1}^{n}t_j(e_j-e_{j+1}).
    \end{split}
    \]
    Case 1: $K \neq J$. The unique positive expression of $e_{K_M}-e_J$ can be written as 
    \[
    e_{K_M}-e_J=\left\{
    \begin{aligned}
        &-\sum\limits_{i=1}^{t}\delta_M(k_i')(e_{k_i'}-e_{k_i'+1})+\sum\limits_{j=1}^{n}t_j(e_j-e_{j+1}), \quad &M\backslash Supp(e_K-e_J)=\emptyset,\\
        &\sum\limits_{l=1}^{n}(e_l-e_{l+1})-\sum\limits_{i=1}^{t}\delta_M(k_i')(e_{k_i'}-e_{k_i'+1})+\sum\limits_{j=1}^{n}t_j(e_j-e_{j+1}), \quad &M\backslash Supp(e_K-e_J) \neq \emptyset.
    \end{aligned}
    \right.
    \]
    so the distance from $e_J$ to $e_{K_M}$ can be expressed by
    \[
    d(e_J,e_{K_M})=\left\{
    \begin{aligned}
        &d(e_J,e_K)-|M|, \quad &M \backslash Supp(e_K-e_J)=\emptyset,\\
        &d(e_J,e_K)-|M|+n, \quad &M\backslash Supp(e_K-e_J). \neq \emptyset
    \end{aligned}
    \right.
    \]
    
   \noindent We noticed that every $M \subset K_{\bullet}$ has a unique partition $M=M' \sqcup M''$ where $M' \subset K_{\bullet} \backslash Supp(e_K-e_J)$ and $M'' \subset K_{\bullet} \cap Supp(e_K-e_J)$. By Lemma \ref{lem2}, $K_{\bullet} \cap Supp(e_K-e_J) \neq \emptyset$, so we consider two subcases as follows\\\\
   Subcase 1: $|K_{\bullet}|=1$ (This is equivalent to that $K$ is an interval). Then $K_{\bullet} \subset Supp(e_K-e_J)$, so the coefficient of $\beta_K$ is 
   \[
   -d(e_J,e_K)+(d(e_J,e_K)-1)=-1.
   \]
   
   \noindent Subcase 2: $|K_{\bullet}| \geq 2$. Then
   \[
   \begin{split}
       \sum\limits_{M \subset K_{\bullet}} (-1)^{1+|M|}d(e_J,e_{K_M}) &=(-1)^{1+|L|}\sum_{\substack{M \subset K_{\bullet}\\M'\neq \emptyset}}\left(d(e_J,e_K)-|M|+n\right)+\sum_{\substack{M \subset K_{\bullet}\\M=M''}}(-1)^{1+|M|}\left(d(e_J,e_K)-|M|\right),\\
       &=\sum\limits_{M \subset K_{\bullet}}(-1)^{1+|M|}(d(e_J,e_K)-|M|+n)-\sum_{\substack{M \subset K_{\bullet}\\M=M''}}(-1)^{1+|M|} \cdot n=0-0=0.
   \end{split}   
   \]
   \noindent Case 2: $K=J$. In this case, 
   \[
   e_{K_M}-e_J=\sum\limits_{l=1}^{n}(e_l-e_{l+1})-\sum\limits_{i=1}^{t}\delta_M(k_i')(e_{k_i'}-e_{k_i'+1})
   \]
   is the unique positive expression of $e_{K_M}-e_J$ for $M \neq \emptyset$, so $d(e_{K_M},e_J)=n-|M|$. Similarly, we consider two subcases.\\\\
   \noindent Subcase 3: $|K_{\bullet}|=1$. The coefficient of $\beta_K$ is $n-1$.\\
   \noindent Subcase 4: $|K_{\bullet}| \geq 2$. The coefficient is 
   \[
       \sum\limits_{M \subset K_{\bullet}} (-1)^{1+|M|}d(e_J,e_{K_M}) =\sum_{\substack{M \subset K_{\bullet}\\L\neq \emptyset}} (-1)^{1+|M|}(n-|M|)=n.
   \]
   
\end{proof}

\begin{lem}
\label{lem2}
    For any $K,J \subset \binom{[n]}{k}$ and $K \neq J$, we have $K_{\bullet} \cap Supp(e_K-e_J) \neq \emptyset$.
\end{lem}
\begin{proof}
    Take the unique positive expression  $e_K-e_J=\sum\limits_{l=1}^{n}t_l(e_l-e_{l+1})$.   Observe that for each $l$, either $|t_l-t_{l+1}|=0$ or $1$. Since $K \neq J$, there exists some $a \in K\backslash J$ such that $t_a \neq 0$. Then we can take the minimal positive integer $b$ such that $t_{a+b}=0$ (Such $b$ do exists since not all the $t_l$ are positive). This implies that $t_{a+b-1}=1$. So the coefficient of $e_{a+b}$ is $-1$ which implies that $a+b \in J\backslash K$. Therefore, there must exist some element of $K_{\bullet}$ such that it lies in the interval $[a, a+b]$.
\end{proof}

\vspace{0.3cm}
\begin{cor}
\label{cor2}
    $(\mathcal{L}_J)_{J \in \binom{[n]}{k}}$ is an $\mathbb{R}$-basis of linear space $\mathfrak{B}_{k,n}=\mathrm{span}_{\mathbb{R}} \left\{\beta_J\;:\;J \in \binom{[n]}{k}\right\}$.
\end{cor}
\begin{proof}
    The transition matrix from $(\beta_J)_{J \in \binom{[n]}{k}}$ to $(n\beta_J-(\beta_{[1,k]}+\beta_{[2,k+1]}+\cdots+\beta_{[n,k-1]}))_{J \in \binom{[n]}{k}}$ is
    \[
    \left(
    \begin{array}{cccc|cccc}
        n &0 &\cdots &0 &0 &0 &\cdots &0\\
        0 &n &\cdots &0 &0 &0 &\cdots &0\\
        \vdots &\vdots &\ddots &\vdots &\vdots &\vdots &\ddots &\vdots \\
        0 &0 &\cdots &n &0 &0 &\cdots &0\\
        \hline
        -1 &-1 &\cdots &-1 &n-1 &-1 &\cdots &-1\\
        -1 &-1 &\cdots &-1 &-1 &n-1 &\cdots &-1\\
        \vdots &\vdots &\ddots &\vdots &\vdots &\vdots &\ddots &\vdots\\
        -1 &-1 &\cdots &-1 &-1 &-1 &\cdots &n-1
    \end{array}
    \right)
    \]
    which is invertible for $n \geq 4$. So $(\mathcal{L}_J)_{J \in \binom{[n]}{k}}$ is a basis.
\end{proof}

\begin{rmk}
    \label{rmk2}
    Similarly to Definition~\ref{def2}, for any $L \subset [n]$ with $0 \leq |L| \leq k-2$, we define $(\mathcal{L}_J^{(L)})_{J \in \binom{[n]\setminus L}{k-|L|}}$ as the linear combination of the family $(\beta_J^{(L)})_{J \in \binom{[n]\setminus L}{k-|L|}}$. Here, the cyclic intervals in the cyclic decomposition of 
$J$ are taken within the base set $[n] \setminus L$, whose cyclic order is inherited from $[n]$ (as explicitly defined in Definition~\ref{def3}). It then follows from Corollary~\ref{cor2} that $(\mathcal{L}_J^{(L)})_{J \in \binom{[n]\setminus L}{k-|L|}}$ forms an $\mathbb{R}$-basis of $\mathfrak{B}_{k,n}^{(L)}$.
\end{rmk}
\vspace{0.5cm}
\begin{ex}
    Take $k=2,\;n=4$ and $J=\{1,3\}$, we have 
    \[
    \begin{split}
        d(e_{13},e_{12})=d(e_{13},e_{34})=1 \quad & d(e_{13},e_{13})=0\\
        d(e_{13},e_{14})=d(e_{13},e_{23})=3 \quad &d(e_{13},e_{24})=2\\
    \end{split}
    \]
    It is not hard to check that $d(e_J,-)$ satisfies the positive tropical $Pl\ddot{u}cker$ relation $d(e_{13},e_{13})+d(e_{13},e_{24})={\rm min}\{d(e_{13},e_{12})+d(e_{13},e_{34}),\;d(e_{13},e_{14})+d(e_{13},e_{23})\}$. The basis $(\mathcal{L}_J)_{J \in \binom{[4]}{2}}$ are
    \[
    \begin{split}
        \mathcal{L}_{12}=-\beta_{12}+\beta_{24} \quad &\mathcal{L}_{14}=-\beta_{14}+\beta_{13} \\
        \mathcal{L}_{23}=-\beta_{23}+\beta_{13}
        \quad &\mathcal{L}_{34}=-\beta_{34}+\beta_{24}\\
        \mathcal{L}_{24}=-\beta_{24}+\beta_{14}
        +\beta_{23}-\beta_{13} \quad &\mathcal{L}_{13}=-\beta_{13}+\beta_{12}
        +\beta_{34}-\beta_{24}
    \end{split}
    \]
    Thus
    \[
    \begin{split}
        \sum\limits_{I \in \binom{[4]}{2}} d(e_J,e_I)\mathcal{L}_I &=(-\beta_{12}+\beta_{24})+3(-\beta_{14}+\beta_{13})+3(-\beta_{23}+\beta_{13})+(-\beta_{24}+\beta_{14}+\beta_{23}-\beta_{13})+(-\beta_{34}+\beta_{24})\\
        &=4\beta_{13}-\beta_{12}-\beta_{23}-\beta_{34}-\beta_{14}.
    \end{split}
    \]
\end{ex}
\vspace{0.5cm}
\begin{defi}
\label{def3}
    Given $L \subseteq [n]$ with $|L|=m\;(0\leq m \leq k-2)$, write $[n]\backslash L=\{l_1,l_2,\cdots,l_{n-m}\}$ with the natural order $1 \leq l_1<l_2<\cdots<l_{n-m} \leq n$. For any $I \subseteq [n-m]$, define $l_I:=\{l_i \;|\;i \in I\}$. We say $l_I$ is \textbf{frozen} in $[n]\backslash L$ if $I$ is a cyclic interval in $[n-m]$; otherwise, it is \textbf{unfrozen}.
\end{defi}
\vspace{0.5cm}
\begin{prop}(\cite{E2})
    \label{prop3}
    For any $j \in [n]\backslash L$,  the boundary map $\partial_j$ acts on the basis $(\mathcal{L}_J^{(L)})_{J \in \binom{[n]\backslash L}{k-|L|}}$ by
    \[
    \partial_j(\mathcal{L}_J^{(L)})=\left\{
    \begin{aligned}
        & \mathcal{L}_{J\backslash\{j\}}^{(L\cup\{j\})} \quad &j \in J\\
        &0 \quad &j \notin J
    \end{aligned}
    \right.
    \]
\end{prop}
\begin{proof}
    To avoid more notations, we suppose that $L=\emptyset$. Let $J^{\bullet}=\{j_1,j_2,\cdots,j_t\}$. Therefore, for any $j\in [n]$, there exists a unique  pair $\{j_a,j_{a+1}\}$ such that $j_a \leq j < j_{a+1}$. Then 
    \[
    \partial_j(\mathcal{L}_J)=\sum\limits_{M \subset J^{\bullet}}(-1)^{1+|M|}\partial_j(\beta_{J_M})
    \]
    If $j \notin J$, we have $\partial_j(\beta_{J_M})=\partial_j(\beta_{J_{M\cup \{j_{a+1}\}}})$ for any $M \subset J^{\bullet} \backslash \{j_{a+1}\}$. Therefore, 
    \[
    \partial_j(\mathcal{L}_J)=\sum\limits_{M \subset J^{\bullet} \backslash \{j_{a+1}\}}((-1)^{1+|M|}+(-1)^{1+|M\cup \{j_{a+1}\}|})\partial_j(\beta_{J_M})=0.
    \]
    Or if $j \in J$, then $\partial_j(\mathcal{L}_J)=\mathcal{L}_{J\backslash\{j\}}^{(\{j\})}$ can be vertified directly.
\end{proof}

\vspace{0.3cm}
\begin{cor}(\cite{E2})
    \label{cor3}
    Given any $\mathbf{c}=(c_I)_{I \in \binom{[n]}{k}} \in \mathbb{R}^{\binom{[n]}{k}}$ and $\mathcal{L}(\mathbf{c}):=\sum\limits_{J \in \binom{[n]}{k}}c_J\mathcal{L}_J$. $\partial_L(\mathcal{L}(\mathbf{c}))$ can be written as a linear combination as 
    \[
    \partial_L(\mathcal{L}(\mathbf{c}))=\sum\limits_{i,j \in [n-k+2]}\pi_{l_il_j}^{(L)}\beta_{l_il_j}^{(L)}
    \]
    where $|L|=k-2$ and $[n] \backslash L=\{l_1,l_2,\cdots,l_{n-k+2}\}$. Choose any unfrozen pair $\{l_i,l_j\}$ (i.e.$|i-j| \geq 2$), the coefficient
    \[
    \pi_{l_il_j}^{(L)}=-(c_{Ll_il_j}-c_{Ll_il_{j+1}}-c_{Ll_{i+1}l_j}+c_{Ll_{i+1}l_{j+1}})
    \]
\end{cor}
\begin{proof}
   The coefficient of $\beta_{l_il_j}^{(L)}$ is from $\mathcal{L}_{l_il_j}^{(L)}$, $\mathcal{L}_{l_{i+1}l_j}^{(L)}$, $\mathcal{L}_{l_il_{j+1}}^{(L)}$ and $\mathcal{L}_{l_{i+1}l_{j+1}}^{(L)}$. Thus by Proposition \ref{prop3}, we get this relation.
\end{proof}

\vspace{0.3cm}

 \noindent Given
\[\beta(\mathbf{c})=\sum\limits_{J \in \binom{[n]}{k}}c_I\beta_I \in \mathfrak{B}_{k,n}
\] 
where $\mathbf{c}=(c_I) \in \mathbb{R}^{\binom{n}{k}}$. Choose some $L \subset [n]$ with $|L|=k-2$, then the image of $\beta(\mathbf{c})$ under the boundary map is
\[
\partial_L(\beta(\mathbf{c}))=\sum\limits_{\{l_i,l_j\} \subset [n]\backslash L } \pi_{l_il_j}^{(L)}\beta_{l_il_j}^{(L)}
\]
where $\pi_{l_il_j}^{(L)}=\sum\limits_{\left\{I \in \binom{[n]}{k}\;:\;\partial_L(\beta_I)=\beta_{l_il_j}^{(L)}\right\}}c_I$. Let the notation $Supp_L\left(\beta(\mathbf{c})\right)$ denote the \textbf{support} of $\beta(\mathbf{c})$ under $\partial_L$ , that is, the set of $\beta_{l_il_j}^{(L)}$ such that $\pi_{l_il_j}^{(L)} \neq 0$.\\

\vspace{0.5cm}
\begin{defi}
\label{def4}
 Denote by $\mathcal{Z}_{k,n}$ the collection of elements $z \in \mathfrak{B}_{k,n}$ that  satisfy
 \begin{enumerate}
     \item  (Compatibility condition) $Supp_L(z)$ induce a positroid subdivision of $\partial_L(\Delta_{k,n})$.
     \item  (Positivity condition) $\pi_{l_il_j}^{(L)} \geq 0$ for any nonfrozen pair $\{l_i,l_j\}$ in $[n] \backslash L$.
 \end{enumerate}
 for any $L \subset [n]$ with $|L|=k-2$.
\end{defi}

\vspace{0.3cm}
The following lemma establishes an explicit connection between $\mathcal{Z}_{k,n}$ and $\rm Trop^+Gr_{k,n}$, linking Early's construction with the positive tropical Grassmannian of Speyer and Williams.

\vspace{0.2cm}
\begin{lem}(\cite{E2})
   \label{lem4}
    Let $\mathbf{c} =(c_I)_{I \in \binom{[n]}{k}} \in \mathbb{R}^{\binom{n}{k}}$. Then $\mathbf{c}$ satisfies the positive tropical Pl\"ucker relation
    \[
    c_{Ll_al_c}+c_{Ll_bl_d}={\rm min}\{c_{Ll_al_b}+c_{Ll_cl_d},c_{Ll_al_d}+c_{Ll_bl_c}\}
    \]
    for any cyclic order $a <b<c<d$ in $[n-k+2]$ if and only if $\mathcal{L}(\mathbf{c}):=\sum\limits_{I}c_I\mathcal{L}_I  \in \mathcal{Z}_{k,n}$. Moreover, this induces a bijection between the positive tropical Grassmannian $\rm Trop^+Gr_{k,n}$ and $\mathcal{Z}_{k,n}$.\\
\end{lem}
\begin{proof}
    Suppose that $\mathbf{c} \in \rm Trop^+Gr_{k,n}$. In particular, we take $a=i,\;b=i+1,\;c=j,\;d=j+1$, then by Corollary \ref{cor3}
    \[
    \pi_{l_il_j}^{(L)}=-(c_{Ll_il_j}-c_{Ll_il_{j+1}}-c_{Ll_{i+1}l_j}+c_{Ll_{i+1}l_{j+1}}) \geq 0
    \]
    for any nonfrozen pair $\{l_i,l_j\}$. It remains to show that $Supp_L(\mathcal{L}(\mathbf{c}))$ induce a positroid subdivision of $\partial(\Delta_{k,n})$. Since $\{l_i,l_j\}$ and $\{l_p,l_q\}$ are not weakly separated  where $ p \in [i+1,j-1]$ and $q \in [j+1,i-1]$, we have
    \[
    \sum_{\substack{p \in [i+1,j-1]\\q \in [j+1,i-1]}} \pi_{l_pl_q}^{(L)}=-(c_{Ll_il_j}-c_{Ll_il_{i+1}}-c_{Ll_jl_{j+1}}+c_{Ll_{i+1}l_{j+1}}).
    \]
    It follows that ${\rm min}\left\{\pi_{l_il_j}^{(L)},\; \sum_{\substack{p \in [i+1,j-1]\\q \in [j+1,i-1]}} \pi_{l_pl_q}^{(L)}\right\}$ is equal to
    \[
    \begin{split}
       &{\rm min}\{-(c_{Ll_il_j}-c_{Ll_il_{j+1}}-c_{Ll_{i+1}l_j}+c_{Ll_{i+1}l_{j+1}}),-(c_{Ll_il_j}-c_{Ll_il_{i+1}}-c_{Ll_jl_{j+1}}+c_{Ll_{i+1}l_{j+1}})\}\\
       =\;&{\rm min} \{(c_{Ll_il_{j+1}}+c_{Ll_{i+1}l_j}),\;(c_{Ll_il_{i+1}}+c_{Ll_jl_{j+1}})\}-(c_{Ll_il_j}+c_{Ll_{i+1}l_{j+1}})\\
       =\;&0.
    \end{split}
    \]
    Thus this proves that $\mathcal{L}(\mathbf{c})$ satisfies (i) in Definition \ref{def4}.\\\\
    \noindent Conversely, suppose that $\mathcal{L}(\mathbf{c}) \in \mathcal{Z}_{k,n}$. Given any $a<b<c<d$ in $[n-k+2]$, we observed that any $\{l_p,l_q\}$ in $A$ and $\{l_{p'},l_{q'}\}$ in $B$ are not weakly separated where
    \[
    \begin{split}
        A:&=\{\{l_p,l_q \}\;|\;p \in [a,b-1],\;q \in [c,d-1]\}\\
        B:&=\{\{l_{p'},l_{q'}\} \;|\;p' \in [b,c-1],\;q' \in [d,a-1]\}
    \end{split}
    \]
    Therefore $\rm min \left\{\sum\limits_{\{l_p,l_q\} \in A}\pi_{l_pl_q}^{(L)},\;\sum\limits_{\{l_{p'},l_{q'}\} \in B}\pi_{l_{p'}l_{q'}}^{(L)} \right\}=0$ by conditions (i) and (ii) in Definition \ref{def4}. As explained on the left of Figure \ref{fig13}, $\pi_{l_pl_q}^{(L)}$ is presented by the labeled square with a number on each corner. The number on the corner $(p,q)$ refers to  the sign of $c_{Ll_pl_q}$ in the expansion of Corollary \ref{cor3}. 
    \begin{figure}[H]
        \centering
        \includegraphics[width=14cm]{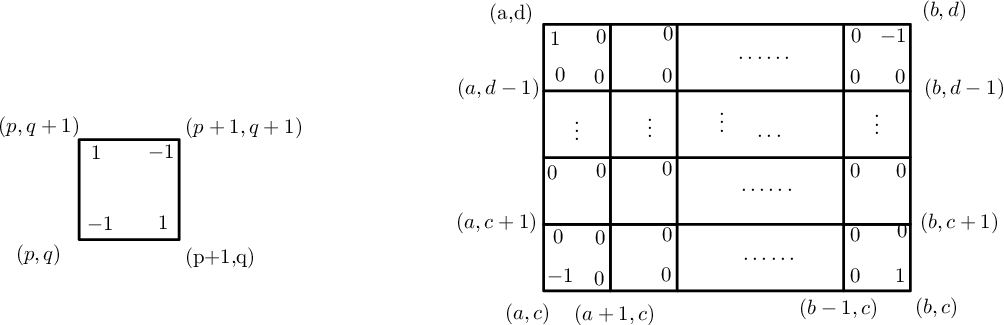}
        \caption{The expansions of $\pi_{l_pl_q}^{(L)}$ and $\sum\limits_{\{l_p,l_q\} \in A}\pi_{l_pl_q}^{(L)}$}
        \label{fig13}
    \end{figure}
    Thus $\sum\limits_{\{l_p,l_q\} \in A}\pi_{l_pl_q}^{(L)}$ can be presented on the right by splicing all the squares together. And we add up all the numbers overlapped at the lattice point. So
    \[
    \begin{split}
        \sum\limits_{\{l_p,l_q\} \in A}\pi_{l_pl_q}^{(L)}=-(c_{Ll_al_c}-c_{Ll_al_d}-c_{Ll_bl_c}+c_{Ll_bl_d})\\
        \sum\limits_{\{l_{p'},l_{q'}\} \in B}\pi_{l_{p'}l_{q'}}^{(L)}=-(c_{Ll_bl_d}-c_{Ll_al_b}-c_{Ll_cl_d}+c_{Ll_al_c})
    \end{split}
    \]
    A straightforward verification shows that $\rm min \left\{\sum\limits_{\{l_p,l_q\} \in A}\pi_{l_pl_q}^{(L)},\;\sum\limits_{\{l_{p'},l_{q'}\} \in B}\pi_{l_{p'}l_{q'}}^{(L)} \right\}=0$ is equivalent to $c_{Ll_al_c}+c_{Ll_bl_d}={\rm min}\{c_{Ll_al_b}+c_{Ll_cl_d},c_{Ll_al_d}+c_{Ll_bl_c}\}$.
\end{proof}

\vspace{0.3cm}
\begin{cor}
    \label{cor4}
    The vector $\mathbf{d}_J:=(d(e_J,e_I))_{I \in \binom{[n]}{k}}$ satisfies the positive tropical  Pl\"ucker relations. Moreover, the subdivison $\mathcal{D}_{\mathbf{d}_J}$ is exactly that induced by translated blade $((1,2,\cdots,n))_{e_J}$.
\end{cor}

\begin{proof}
    Since $h(x-e_J)={\rm min}\{L_1(x-e_J), L_2(x-e_J), \cdots, L_n(x-e_J)\}$ is concave on the hyperplane $\sum\limits_{i=1}^{n}x_i=k$, we know that $-h(x-e_J)$ is convex. By Proposition \ref{prop5}, the intersections of all the affine areas of $-h(x-e_J)$ with $\Delta_{k,n}$ form the subdivision induced by $((1,2,\cdots,n))_{e_J}$. On the other hand, because of the convexity of $-h(x-e_J)$ and $-h(e_J,e_I)=d(e_J,e_I)$, the projection of any lower face of $\rm convex(\{(e_I,d(e_J,e_I))\})$ is exactly the intersection of some affine area of $-h(x-e_J)$ with $\Delta_{k,n}$.
\end{proof}

\vspace{0.3cm}
\begin{rmk}
    \label{3}
    From Corollary \ref{cor4} together with the equation in Proposition \ref{prop4}
    \[
    \sum\limits_{I \in \binom{[n]}{k}} d(e_J,e_I)\mathcal{L}_I=n\beta_J-(\beta_{[1,k]}+\beta_{[2,k+1]}+\cdots+\beta_{[n,k-1]}),
    \]
    $\beta_J$ can be regarded as the positroid subdivision induced by $((1,2,\cdots,n))_{e_J}$. Also, we ignore the negative terms on the right since $\beta_I$ induces a trivial subdivision of $\Delta$ when $I$ is a cyclic interval. Through these observations, it is not hard to see that $\mathcal{Z}_{k,n}$ contains a linearity spcae $\mathcal{F}_{k,n}:=\rm span \{\beta_I\;|\;I\;is\;frozen\}$.
\end{rmk}

\vspace{0.3cm}
\begin{defi}
    Let $\overline{\mathcal{Z}_{k,n}}:=\mathcal{Z}_{k,n}/ \mathcal{F}_{k,n}$ be the quotient fan modulo the linearity space. Then evey element $x \in \overline{\mathcal{Z}_{k,n}}$ can be uniquely written as a weighted blade arrangement of nonfrozen translated blades.
\end{defi}

\vspace{1cm}

\subsection{The first theorem}\quad\\
\label{subsec1}
From now on, we restrict our scope to weighted blade arrangements in $\overline{\mathcal{Z}_{k,n}}$, by which we mean that all terms of frozen translated blades are omitted. Let us first consider the case $k=2$.
\vspace{0.2cm}
\begin{prop}
    \label{prop6}
    All the maximal cones of $\overline{\mathcal{Z}_{2,n}}$ are simplicial with dimension $n-3$. And they are of the form $\mathcal{P}_W:=\{\sum_{J \in W}\mathbb{R}_{\geq 0}\beta_J\}$ for some nonfrozen maximal w-collection $W$.
\end{prop}

\begin{proof}
    Any finest positroidal subdivision of$\Delta_{2,n}$ contains $n-2$ top dimensional cells (see \cite{SW2})  and $n-3$ facets which split these cells (2-splits). Every translated blade to a non-frozen vertex induces a 2-split. Therefore we have 
weakly separated collection of $n-3$ vertices corresponding to these splits. But $n-3$ is a maximal cardinality of a weakly separated collection of two-elements sets. This implies that the cone spanned by  a such collection of weighted blades forms a simplicial cone of
$\overline{\mathcal{Z}_{2,n}}$.  Moreover, there is a bijection between maximal weakly separated sets of two elements sets and vertices of the associahedron (see Fomin and Zelevinsky \cite{FZ} ) .  Therefore all simplicial cones of $\overline{\mathcal{Z}_{2,n}}$ are of such a form.
\end{proof}

\vspace{0.3cm}

From Proposition \ref{prop6}, we know that every maximal cone in $\overline{\mathcal{Z}_{2,n}}$ is simplicial. However when $k \geq 3$, it is not the case. Here we provide a necessary and sufficient condition for $\mathcal{P}_{W}:=\{\sum_{J \in W}\mathbb{R}_{\geq}\beta_J\}$ being a maximal simplicial cone in $\overline{\mathcal{Z}_{k,n}}$ where $W$ is a nonfrozen maximal $w$-collection.

Given a collection of $k$-sets $X \subseteq \binom{[n]}{k}$ and $L  \in \binom{[n]}{m}$ ($0 \leq |L| \leq k-2$), let
\[
\partial_L(X):=\{\partial_L(J)\;|\;J \in X\}.
\]
Define
\[
\mathscr{S}:=\{X \subseteq \binom{[n]}{k} \;|\;\partial_L(X) \text{ contains a nonfrozen maximal w-collection in} \; \binom{[n]\backslash L}{2},\; \forall L \in \binom{[n]}{k-2}\}
\]
as a poset of nonfrozen $w$-collections ordered by inclusion.
\begin{thm}
\label{thm1}
Let $W$ be a nonfrozen maximal $w$-collection. The following statements are equivalent:
\begin{enumerate}
    \item $W$ is minimal in the poset $\mathscr{S}$,
    \item $\mathcal{P}_{W}$ is a maximal simplicial cone in the fan $\overline{\mathcal{Z}_{k,n}}$.
\end{enumerate}
\end{thm}

\begin{proof}
The equivalence between (1) and (2) is directly obtained by Proposition \ref{prop6} when $k=2$. So we assume $k \geq 3$.

(1) $\Rightarrow$ (2).  Here exists a maximal cone $\mathcal{M}$ of $\overline{\mathcal{Z}_{k,n}}$ such that $\mathcal{P}_{W} \subseteq \mathcal{M}$ and all generators of $\mathcal{P}_{W}$ are generators of $\mathcal{M}$. We show that $\mathcal{P}_{W} = \mathcal{M}$, i.e., every proper face $\mathcal{P}_{W'}$ ($W' \subsetneqq W$) is disjoint from the interior of $\mathcal{M}$.

This follows from the minimality of $W$ in $\mathscr{S}$, combined with Proposition \ref{prop6}, which ensures the existence of $L_0 \in \binom{[n]}{k-2}$ such that $\{\beta_{J}^{(L_0)} \mid J \in \partial_{L_0}(W')\}$ does not induce a finest positroid subdivision of $\partial_{L_0}(\Delta_{k,n}) = \Delta_{2,[n]\backslash L_0}$. That is, there exists an octahedron in $\partial_{L_0}(\Delta_{k,n})$ that remains unsubdivided.

By Theorem \ref{thm3} and the secondary fan structure of $\overline{\mathcal{Z}_{k,n}}$, we conclude that $\{\beta_J \mid J \in \mathcal{P}_{W'}\}$ does not induce a finest positroid subdivision of $\Delta_{k,n}$. Hence, $\mathcal{P}_{W'}$ is disjoint from the interior of $\mathcal{M}$,i.e., $\mathcal{P}_{W} = \mathcal{M}$.

(2) $\Rightarrow$ (1). When $\mathcal{P}_W$ is maximal and simplicial, any $\mathcal{P}_{W'}$ ($W' \subsetneqq W$) is a proper face of $\mathcal{P}_W$. By the secondary fan structure of $\overline{\mathcal{Z}_{k,n}}$, there always exists an octahedron $\{L_0ab, L_0bc, L_0cd, L_0ad, L_0ac, L_0bd\}$ that is not subdivided by the refinement induced by $\{\beta_J \mid J \in \mathcal{P}_{W'}\}$.

Therefore, according to Proposition \ref{prop6}, $\partial_{L_0}(W')$ does not contain a nonfrozen maximal w-collection, i.e., $W' \notin \mathscr{S}$. This implies that $W$ is minimal in $\overline{\mathcal{Z}_{k,n}}$.
\end{proof}

\vspace{0.3cm}

\begin{rmk}
    Based on Theorem \ref{thm2} in the subsequent section, it can be established that any nonfrozen maximal $w$-collection $W$ is contained within $\mathscr{S}$, although not every such $W$ constitutes a minimal element in $\mathscr{S}$. Consequently, the enumeration of all maximal simplicial cones of this type remains an open problem.
\end{rmk}

\vspace{0.2cm}
\begin{ex}
    When $k=3$ and $n=6$, we consider the nonfrozen w-collection $W=\{124,125,134,145\}$. Take the boundary map to secondary hypersimplex faces of $\Delta_{3,6}$,
    \begin{equation*}
        \begin{split}
            \partial(\beta_{124}) &=\beta_{24}^{(1)}+\beta_{14}^{(2)}+\beta_{24}^{(5)}+\beta_{24}^{(6)},\\
            \partial(\beta_{125}) &=\beta_{25}^{(1)}+\beta_{15}^{(2)}+\beta_{25}^{(6)},\\
            \partial(\beta_{134})&=\beta_{14}^{(2)}+\beta_{14}^{(3)}+\beta_{13}^{(4)},\\
            \partial(\beta_{145})&=\beta_{15}^{(2)}+\beta_{15}^{(3)}+\beta_{15}^{(4)}+\beta_{14}^{(5)}.
        \end{split}
    \end{equation*}
It can be checked that $W$ is minimal in $\mathscr{S}$, thus $\mathcal{P}_{W}$ comprise  a maximal simplicial cone in $\overline{\mathcal{Z}_{3,6}}$.
\end{ex}
Below is an example serving to illustrate that the positive weighted blade arrangements arising from certain nonfrozen maximal $w$-collections do not necessarily yield a maximal simplicial cone within $\overline{\mathcal{Z}_{k,n}}$.

\vspace{0.3cm}
\begin{ex}
\label{exam1}
    Still in $\overline{\mathcal{Z}_{k,n}}$, consider a nonfrozen maximal $w$-collection  $W=\{135,235,145,136\}$. Then
    \begin{equation*}
        \begin{split}
            \partial(\beta_{135}) &=\beta_{35}^{(1)}+\beta_{15}^{(2)}+\beta_{15}^{(3)}+\beta_{13}^{(4)}+\beta_{13}^{(5)}+\beta_{35}^{(6)},\\
            \partial(\beta_{235}) &=\beta_{35}^{(1)}+\beta_{35}^{(2)}+\beta_{25}^{(3)}+\beta_{35}^{(6)},\\
            \partial(\beta_{145})&=\beta_{15}^{(2)}+\beta_{15}^{(3)}+\beta_{15}^{(4)}+\beta_{14}^{(5)},\\
            \partial(\beta_{136})&=\beta_{36}^{(1)}+\beta_{13}^{(4)}+\beta_{13}^{(5)}+\beta_{13}^{(6)}.
        \end{split}
    \end{equation*}
        It follows directly from the above computation that $W$ contains a proper subset $W'=\{235,145,136\}$ which lies in $\mathscr{S}$. Thus $W$ is not minimal in $\mathscr{S}$, i.e., $\mathcal{P}_{W}$ is not a maximal simplicial cone of $\overline{\mathcal{Z}_{k,n}}$.
    
    And here we explicitly write the generators of this cone in standard coordinate $\{\mathcal{L}_I\}$ of $\mathbb{R}^{\binom{6}{3}}$ where $I \in \binom{[6]}{3}$.\\
     \begin{equation*}
         \begin{split}
             \mathbb{R}_{\geq 0}\beta_{135} =&\mathbb{R}_{\geq 0}(3\mathcal{L}_{123}+2\mathcal{L}_{124}+\mathcal{L}_{125}+6\mathcal{L}_{126}+\mathcal{L}_{134}+5\mathcal{L}_{136}+5\mathcal{L}_{145}+4\mathcal{L}_{146}+3\mathcal{L}_{156}+\\
             &6\mathcal{L}_{234}+5\mathcal{L}_{235}+4\mathcal{L}_{236}+4\mathcal{L}_{245}+3\mathcal{L}_{246}+2\mathcal{L}_{256}+3\mathcal{L}_{345}+2\mathcal{L}_{346}+\mathcal{L}_{356}+6\mathcal{L}_{456})\\\\
             \mathbb{R}_{\geq 0}\beta_{235} =&\mathbb{R}_{\geq 0}(4\mathcal{L}_{123}+3\mathcal{L}_{124}+2\mathcal{L}_{125}+7\mathcal{L}_{126}+2\mathcal{L}_{134}+1\mathcal{L}_{135}+6\mathcal{L}_{136}+6\mathcal{L}_{145}+5\mathcal{L}_{146}+\\
             &6\mathcal{L}_{156}+\mathcal{L}_{234}+5\mathcal{L}_{236}+5\mathcal{L}_{245}+4\mathcal{L}_{246}+3\mathcal{L}_{256}+5\mathcal{L}_{345}+3\mathcal{L}_{346}+2\mathcal{L}_{356}+7\mathcal{L}_{456})\\\\
             \mathbb{R}_{\geq 0}\beta_{145} =&\mathbb{R}_{\geq 0}(4\mathcal{L}_{123}+3\mathcal{L}_{124}+2\mathcal{L}_{125}+7\mathcal{L}_{126}+2\mathcal{L}_{134}+\mathcal{L}_{135}+6\mathcal{L}_{136}+5\mathcal{L}_{146}+4\mathcal{L}_{156}+\\
             &7\mathcal{L}_{234}+6\mathcal{L}_{235}+5\mathcal{L}_{236}+5\mathcal{L}_{245}+4\mathcal{L}_{246}+3\mathcal{L}_{256}+4\mathcal{L}_{345}+3\mathcal{L}_{346}+2\mathcal{L}_{356}+7\mathcal{L}_{456})\\\\
             \mathbb{R}_{\geq 0}\beta_{136} =&\mathbb{R}_{\geq 0}(4\mathcal{L}_{123}+3\mathcal{L}_{124}+2\mathcal{L}_{125}+1\mathcal{L}_{126}+2\mathcal{L}_{134}+\mathcal{L}_{135}+6\mathcal{L}_{145}+5\mathcal{L}_{146}+4\mathcal{L}_{156}+\\
             &7\mathcal{L}_{234}+6\mathcal{L}_{235}+5\mathcal{L}_{236}+5\mathcal{L}_{245}+4\mathcal{L}_{246}+4\mathcal{L}_{256}+4\mathcal{L}_{345}+3\mathcal{L}_{346}+2\mathcal{L}_{356}+7\mathcal{L}_{456})\\\\
         \end{split}
     \end{equation*}
    The above equations hold modulo the linearity space $\rm span_{\mathbb{R}}\{\sum\limits_{i \in I} \mathcal{L}_I\;|\;i=1,\cdots,6\}$ of $\rm Trop^+Gr_{3,6}$.
\end{ex}

\vspace{0.2cm}

\begin{rmk}
    To further illustrate the connection with Speyer and Williams \cite{SW}, we recover their polymake-based computational results within the framework of weighted blades.

  The fan of $\operatorname{Trop}^+ Gr_{3,6}$ is a coarsening of the normal fan of the type $D_4$ associahedron, with f-vectors $f=(f_1,f_2,f_3,f_4)=(16,66,98,48)$ ($f_i$ denotes the number of $i$-dimensional faces) and $f=(16,66,100,50)$, respectively. The coarsening is induced by two non-simplicial bipyramid maximal cones; subdividing these bipyramid into tetrahedra recovers the finer type $D_4$ fan.

    Consider two nonfrozen maximal $w$-collections $W_1=\{135,235,145,136\}$ and $W_2=\{246,346,256,124\}$. By Theorem \ref{thm2}, both of them are not minimal in $\mathscr{S}$, i.e. $\mathcal{P}_{W_1}$  and $\mathcal{P}_{W_2}$ are not maximal simplicial cones of $\operatorname{Trop}^+ Gr_{3,6}$. Therefore, the maximal cones $\mathcal{M}_1$ and $\mathcal{M}_2$ (see Figure \ref{fig22}), which contain $\mathcal{P}_{W_1}$ and $\mathcal{P}_{W_2}$ respectively, also require one additional generator each, namely $-\beta_{135}+\beta_{145}+\beta_{136}+\beta_{235}$ and $-\beta_{246}+\beta_{256}+\beta_{124}+\beta_{346}$.
    \begin{figure}[H]
        \centering
        \includegraphics[width=10cm]{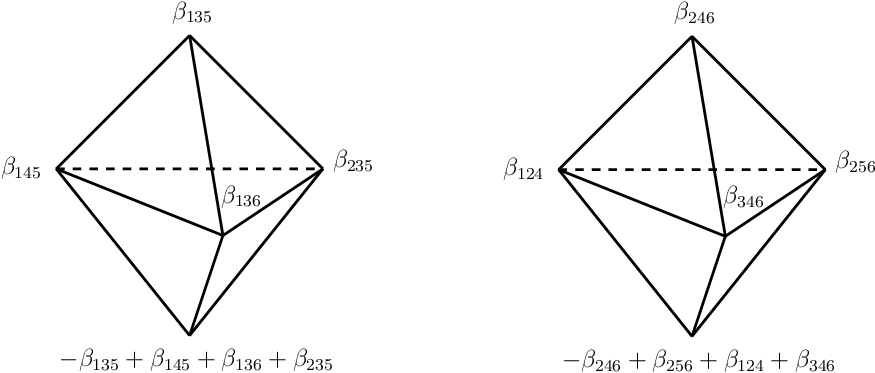}
        \caption{Two non-simplicial bipyramid maximal cones in $\overline{\mathcal{Z}_{3,6}}$.}
        \label{fig22}
    \end{figure}
    To see why \( -\beta_{135} + \beta_{145} + \beta_{136} + \beta_{235} \) is a generator of $\overline{\mathcal{Z}_{3,6}}$ , recall from Remark \ref{rmk1} that weighted blades can be viewed as weighted sums of indicator functions. We represent \( \Delta_{3,6} \) by a circle, where its codimension-1 and codimension-2 faces are depicted as line segments and points, respectively. The numbers 
    on these faces indicate the values of the weighted blade. It is then evident from Figure \ref{fig23} that 
\( -\beta_{135} + \beta_{145} + \beta_{136} + \beta_{235} \) induces a 3-split.
\begin{figure}[H]
        \centering
        \includegraphics[width=13cm]{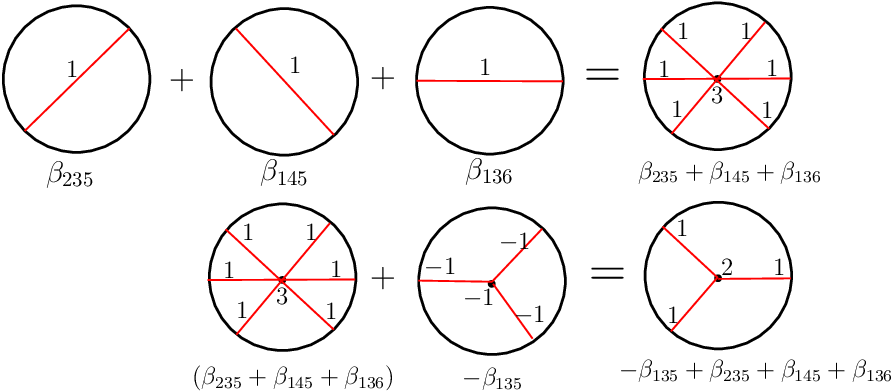}
        \caption{The weighted blade $-\beta_{135} + \beta_{145} + \beta_{136} + \beta_{235}$}
        \label{fig23}
    \end{figure}
\end{rmk}

\vspace{0.2cm}

\vspace{1cm}
\section{Boundary maps on reduced plabic graphs and the second theorem}
\label{sec3}
This section introduces basic properties of reduced plabic graphs from \cite{OPS}. We then recast the boundary map $\partial_j$ on a maximal $w$-collection $W$ as a series of manipulations on its associated reduced plabic graph $\Sigma_0(W)$. To prove Theorem \ref{thm2}, we show that the face labels of $\partial_j(\Sigma_0(W))$ precisely yield $\partial_j(W)$ and that $\partial_j(\Sigma_0(W))$ is again a bipartite reduced plabic graph. This result implies that the translated blades from any maximal w-collection induce a finest positroid subdivision of $\Delta_{k,n}$. We conclude by characterizing flips between these finest subdivisions via pairs of adjacent maximal cones in $\overline{\mathcal{Z}_{k,n}}$.

\begin{thm}
\label{thm2}
    Let $W$ be a maximal w-collection in $\binom{[n]}{k}$, then $\partial_j(W)$ is also a maximal w-collection in $\binom{[n]\backslash \{j\}}{k-1}$.
\end{thm}

The proof of this Theorem will be conducted in the end of this section by using plabic graphs.

\vspace{0.3cm}

\subsection{The boundary maps on plabic graphs}

\begin{defi}
    A plabic graph  is a planar graph $G$ embedded in a closed disk $\mathbf{D}$ satisfying the following conditions
    \begin{enumerate}
        \item No edges cross each other.
        \item Each internal vertex is colored black or white.
        \item Each internal vertex is connected by a path to some boundary vertex.
        \item The boundary vertices are labeled by $1,2,\cdots,n$ in clockwise order for some $n \in \mathbb{Z}_{>0}$.
        \item Each boundary vertex is incident to exactly one internal vertex by an edge.
    \end{enumerate}
\end{defi}

\vspace{0.3cm}
\begin{defi}
    Two plabic graphs $G$ and $G'$ are said to be  move equivalent to each other if $G$ and $G'$ can be related to each other via a sequence of the following local moves (M1), (M2) and (M3):\\\\
    (M1) (The square move) Change the colors of vertices of a quadrilateral face whose vertices are trivalent as Figure \ref{fig10}.
    \begin{figure}[H]
        \centering
        \includegraphics[width=6cm]{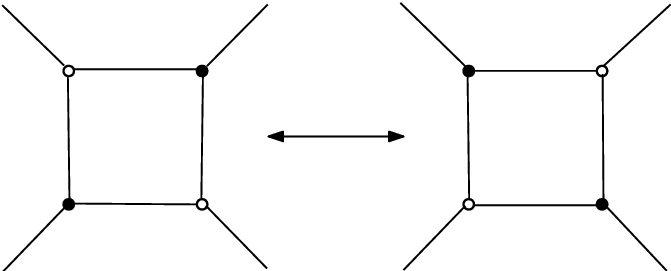}
        \caption{}
        \label{fig10}
    \end{figure}
   \noindent (M2) Contract an edge of connecting two internal vertices of the same color or split an internal vertex into two vertices of the same color joined by an edge as Figure \ref{fig11}
    \begin{figure}[H]
        \centering
        \includegraphics[width=10cm]{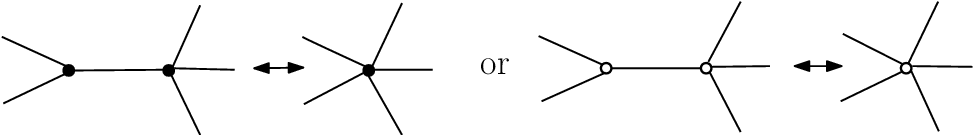}
        \caption{}
        \label{fig11}
    \end{figure}
    \noindent(M3) Remove a bivalent vertex and merge the edge adjacent to it or insert a bivalent in the middle of an edge as Figure \ref{fig12}.
    \begin{figure}[H]
        \centering
        \includegraphics[width=8cm]{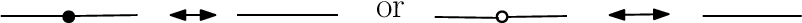}
        \caption{}
        \label{fig12}
    \end{figure}
\end{defi}

In \cite{OPS}, Oh et al. utilized a technique called \textbf{plabic tilings}. Given a maximal $w$-collection $W$, they constructed the corresponding plabic tiling and took its dual graph to obtain a bipartite reduced plabic graph $\Sigma_0(W)$. Conversely, the target labels on all faces of $\Sigma_0(W)$ exactly recover the collection $W$. That is, $W$ and $\Sigma_0(W)$ are in a one-to-one correspondence.	
Here, we will not delve into the detailed construction of plabic tilings or the specific rules for target labels; for details, see \cite{OPS}.

\vspace{0.3cm}
\begin{defi}
    Let $W$ be a maximal $w$-collection in $\binom{[n]}{k}$, and $\Sigma_0(W)$ be the bipartite reduced plabic graph obtained from the duality of plabic tiling. Define $\Sigma(W)$ to be the equivalence class of $\Sigma_0(W)$ under moves (M2) such that every vertex has degree at least 3. (Obviously $\Sigma_0(W) \in \Sigma(W)$ by the definition of plabic tiling)
\end{defi}


\vspace{0.3cm}
Within the framework of plabic tiling construction in \cite{OPS} or the combined tiling in \cite{DKK2}, every black vertex in $\Sigma_0(W)$ corresponds to a $(k+1)$-set $\mathcal{B}$. The labels of the faces surrounding $\mathcal{B}$ in clockwise order are precisely the $k$-sets $\{ \mathcal{B} \setminus b_1, \mathcal{B} \setminus b_2, \dots, \mathcal{B} \setminus b_s \}$ in $W$, where $b_1 < b_2 < \cdots < b_s$ follows the cyclic order. Similarly, every white vertex corresponds to a $(k-1)$-set $\mathcal{W}$, with its adjacent face labels given by $\{ \mathcal{W} \cup w_1, \mathcal{W} \cup w_2, \dots, \mathcal{W} \cup w_r \}$ in clockwise order, where $w_1 < w_2<\cdots <w_r$.

\vspace{0.3cm}
\begin{defi}
    For any $i \in [n]$, $<_i$ denotes a linear order on $[n]$ as follows:
    \[
    i <_i i+1 <_i i+1 \cdots <_i i-1
    \]
    and for any subset $S \subset [n]$, we use the notation $M_l^i(S)$ to denote the set of minimal $l$ elements in $S$ under $<_i$.\\
\end{defi}

\begin{defi}
\label{def1}
    Without loss of generality, we take $j=n$. The plabic graph $\partial_n(\Sigma_0(W))$ is obtained from $\Sigma_0(W)$ through the following steps
    \begin{enumerate}
        \item Delete the vertex labeled $n-k$ on the boundary and the unique edge adjacent to it, but preserve the internal vertex of this edge.
        \item Delete all the edge with vertices $\mathcal{B}$ and $\mathcal{W}$ such that $\mathcal{B} \backslash \mathcal{W}=M_2^n(\mathcal{B})$.
        \item Relabel the boundary vertices $n,1,2,\cdots,n-k-1$ by $1,2,\cdots,n-k$ clockwise.
        \item Delete all the single points and bivalent points, use moves (M2) to get a bipartite plabic graph.
    \end{enumerate}
\end{defi}
We use the same notation as the Appendix of \cite{DKK3}. The collection $W_0=\mathcal{I}_n^k \cup \mathcal{S}_n^k$, where $\mathcal{I}_n^k$ consists of the intervals of size $k$ and $\mathcal{S}_n^k$ consists of the sets of size $k$ represented as the union of two nonempty intervals $[1,m] \cup [l,l+k-m-1]$ with $l \geq m+2$ and $l+k-m-1 \leq n-1$. To prove Theorem \ref{thm2}, we first prove a special case then we show that the maximality of $\partial_j(W_0)$ preserves under flips.\\

\vspace{0.3cm}
\subsection{Some lemmas for preparations}
\begin{lem}
\label{lem6}
    Let $W_0=\mathcal{I}_n^k \cup \mathcal{S}_n^k$, then $\partial_n(\Sigma_0(W_0)) \cong \Sigma_0(\partial_n(W_0))$
\end{lem}

\begin{proof}
    
    We observe that $\partial_n(W_0)$ consists of cyclic intervals of size $k-1$ in $[n-1]$, along with sets of size $k-1$ that can be expressed in the form $[2, m] \cup [l, l + k - m - 1]$. It is clear that $\partial_n(W_0)$ forms a maximal $w$-collection in $\binom{[n-1]}{k-1}$.
    
    Now consider the local area of each face $F$ in $\partial_n(\Sigma_0(W_0))$, where the \textbf{local areas} denoted by $L(F)$ consists of $F$ together with all faces sharing a vertex or an edge with $F$. We classify these local areas $L(F)$ into four cases:\\

    \noindent Case 1: As in Figure \ref{fig1}, the left is the local area of face $F=\{1\} \cup [l,k+l-2]$ with $l\geq 4$.  Since $\partial_n([l-1,k+l-2])=\partial_n(\{1\}\cup [l,k+l-2])=[l,k+l-2]$ and $\partial_n([l,k+l-1])=\partial_n(\{1\}\cup [l+1,k+l-1])=[l+1,k+l-1]$, We remove the two edges that separate these two pairs of faces, respectively. Then, applying the remaining steps in Definition \ref{def1}, we obtain the graph shown on the right.
    \begin{figure}[H]
        \centering
        \includegraphics[width=15cm]{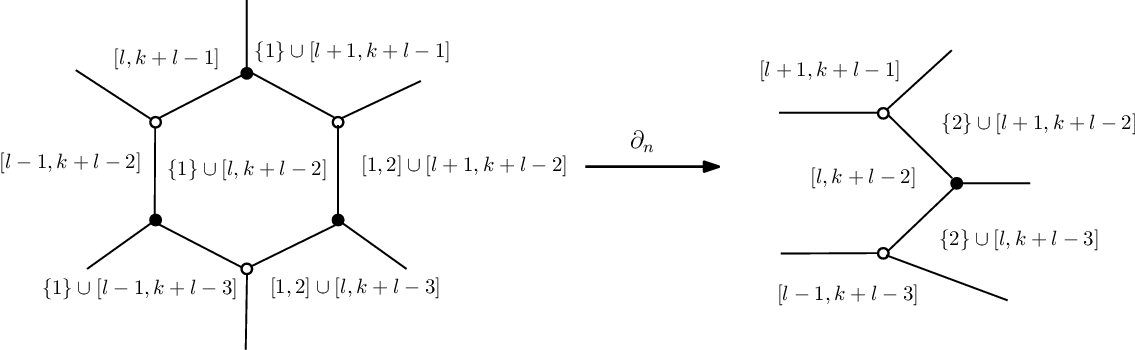}
        \caption{ Boundary map on local area of $F=\{1\} \cup [l,k+l-2]$ with $l \geq 4$}
        \label{fig1}
    \end{figure}

    \noindent Case 2: As in Figure \ref{fig2}, the left is the local area of face $F=\{1\} \cup [3,k+1]$. Since $\partial_n([2,k+l])=\partial_n(\{1\}\cup [3,k+l])=[3,k+l]$ and $\partial_n([3,k+2])=\partial_n(\{1\}\cup [4,k+2])=[4,k+2]$, the same as Case 1, we delete these edges to combine the adjacent faces into one face.
    \begin{figure}[H]
        \centering
        \includegraphics[width=15cm]{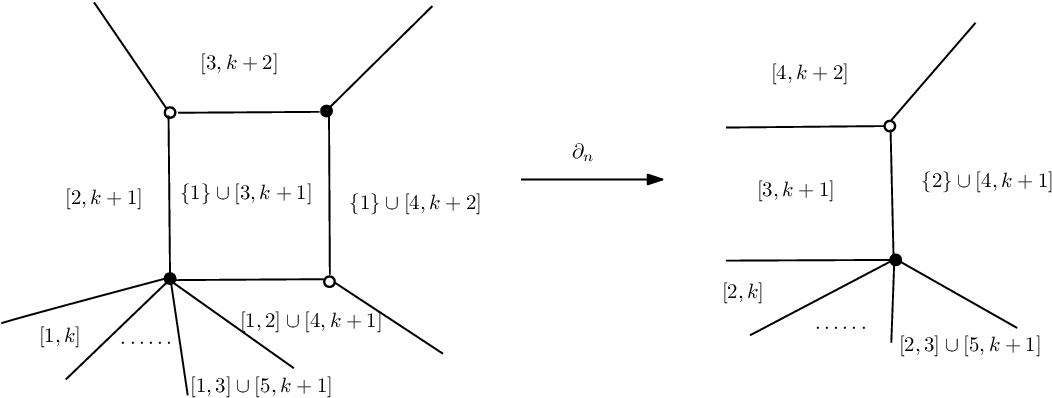}
        \caption{ Boundary map on local area of $F=\{1\} \cup [3,k+l]$}
        \label{fig2}
    \end{figure}
    \noindent Case 3: As shown in Figure \ref{fig3} (left), the local area corresponds to face $F = [1, m] \cup [l, k + l - m - 1]$, where $m \geq 2$ and $l \geq m + 3$. By definition, no adjacent faces within $L(F)$ will be merged into a single face. Therefore, it is only necessary to relabel each face while leaving the graph structure unchanged.
    \begin{figure}[H]
        \centering
        \includegraphics[width=15cm]{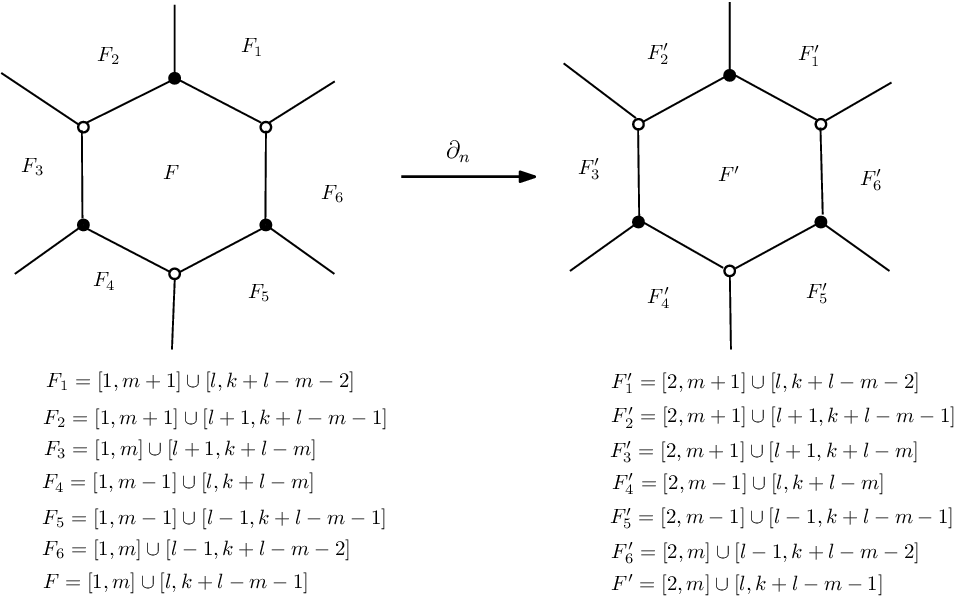}
        \caption{Boundary map on local area of $F=[1,m] \cup [l,k+l-m-1]$ with $m \geq 2$ and $l \geq m+3$}
        \label{fig3}
    \end{figure}
\noindent Case 4: As in Figure \ref{fig4}, the left is the local area of face $F=[1,m] \cup [m+2,k+1]$ with $m \geq 3$. The only edge that need to be deleted is the edge that separates face $F_1=\{1\} \cup [3,k+1][2,k+1]$ and face $F_2=[2,k+1]$. So $F_1$ and $F_2$ are united into $F'_{1,2}$ on the right.
\begin{figure}[H]
    \centering
    \includegraphics[width=16cm]{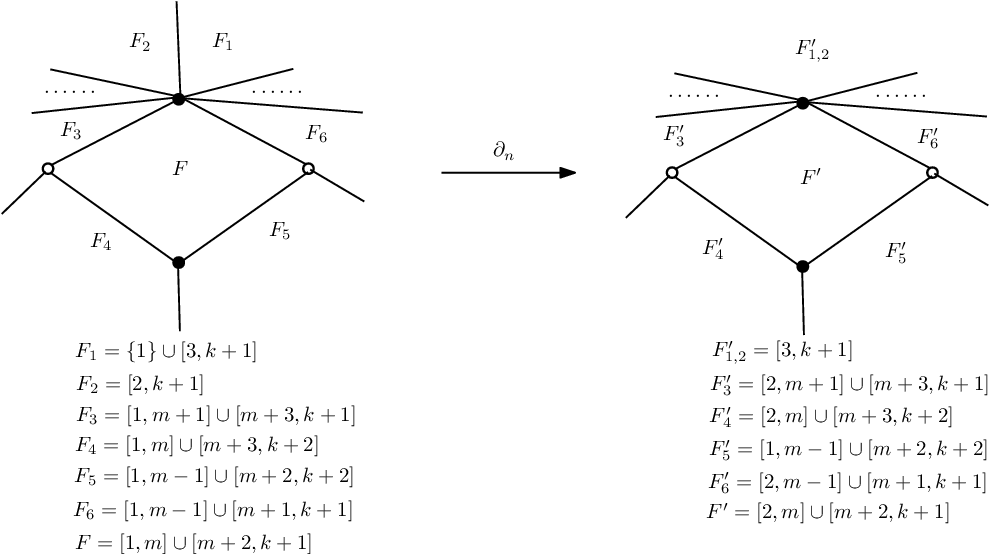}
    \caption{Boundary map on local area of $F=[1,m] \cup [m+2,k+l]$ with $m \geq 3$}
    \label{fig4}
\end{figure}
Given any element $F' \in \partial_n(W_0)$, we observed that the local area of $F'$ in $\Sigma_0(\partial_n(W_0))$ is exactly presented on the right of some case above. Since a plabic graph is a planar graph, the connection of these local areas must be unique. So $\partial_n(\Sigma_0(W_0)) \cong \Sigma_0(\partial_n(W_0))$.
\end{proof}

\begin{figure}[H]
    \centering
    \includegraphics[width=12cm]{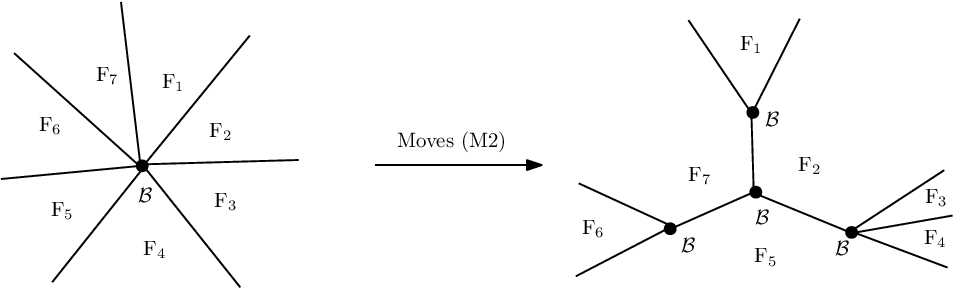}
    \caption{Labels of black points in $G$}
    \label{fig5}
\end{figure}

\begin{rmk}
We can extend the operations in Definition \ref{def1} to the equivalence class $\Sigma(W)$. Since every $G \in \Sigma(W)$ is obtained from $\Sigma_0(W)$ by a series of moves (M2), a black (white) point in $\Sigma_0(W)$ may be extended to a tree with black (white) vertives. We can assign the same label on every point of this tree as the point in $\Sigma_0(W)$. So we can still do the operations in Definition \ref{def1}. This process can be illustrated in Figure \ref{fig5}.
\end{rmk}

For an extension Lemma \ref{lem6} from $W_0$ to any maximal $w$-collection in $[n]\choose k$, we need the fact that any two maximal $w$-collection can be joined by a sequence of flips (cluster mutations). Recall  (Leclerc-Zelevinsky) that if a $w$-collection
$W$ contain five  sets $Lab$, $Lcd$, $Lac$,  $Lad$, $Lbc$, then the collection $W'=W\setminus \{Lac\} \cup \{Lbd\}$ is also weakly separated. The transformation $W\to W'$ is called a \textbf{flip}.
\begin{thm}(\cite{DKK},\cite{OPS})
    Any two maximal w-collection in  $[n]\choose k$ can be joined by a sequence of flips.
\end{thm}

\begin{figure}[H]
    \centering
    \includegraphics[width=5cm]{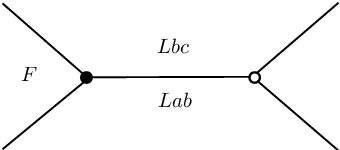}
    \caption{A forbidden structure in plabic graph}
    \label{fig8}
\end{figure}

\vspace{0.2cm}
\begin{lem}
\label{lem1}
    For any $G \in \Sigma(W)$, the following forbidden structure (see Figure \ref{fig8}) with $\mathcal{B} \backslash \mathcal{W}=M_2^n(\mathcal{B})$ will not appear in $G$ ($a,b,c$ are cyclic ordered in $[n] \backslash L$) .
\end{lem}
\begin{proof}
    Suppose that this structure appear in $G$, since $\mathcal{B} \backslash \mathcal{W}=\{Labc\}\backslash \{Lb\}=\{a,c\}=M_2^n(\mathcal{B})$, we have either $(a,c) \cap Lb=\emptyset$ or $(c,a) \cap Lb=\emptyset$. But $a,b,c$ are cyclic ordered, so $(a,c) \cap Lb \neq \emptyset$. However, since evrery internal point in $G$ has degree no less than three, there exists a face on the left of $\mathcal{B}$ such that $F=\{Labc\}\backslash \{x\}$ where $x \in (c,a)$. There is a contradiction to the fact $(c,a) \cap Lb=\emptyset$.
\end{proof}

\vspace{0.3cm}
\begin{lem}
    \label{lem7}
    Let $W$ be any maximal w-collection, then $\partial_n(\Sigma_0(W)) \cong \Sigma_0(\partial_n(W))$. Moreover, the induced subgraph (1-dimensional subcomplex) of points $\partial_n^{-1}|_{W}(\partial_n(J)):=\{I \in W\;|\;\partial_n(I)=\partial_n(J),\; J \in W\}$ is connected in the plabic tiling of $W$.
\end{lem}
\begin{proof}
    We start with $W_0=\mathcal{I}_n^k \cup \mathcal{S}_n^k$, we have known that $\partial_n(\Sigma_0(W_0)) \cong \Sigma_0(\partial_n(W_0))$. We will show that this relation holds under flips. Let $W$ be a maximal w-collection satisfying $\partial_n(\Sigma_0(W)) \cong \Sigma_0(\partial_n(W))$ and $W'=(W\backslash\{Lac\})\cup \{Lbd\}$. We noticed that, for any fixed black point $\mathcal{B}$ in $G \in \Sigma(W)$, there are at most one edge adjacent to $\mathcal{B}$ can be deleted under the boundary map $\partial_n$. Thus we only need to consider two cases up to symmetry. These two cases depend on how many edges are deleted in the local area formed by $Lab,Lbc,Lcd,Lad,Lac (Lbd)$. \\\\
    Case 1: As in figure \ref{fig6}, only one edge is deleted in step 2 in the local area formed by $Lab,Lbc,Lcd,Lad,Lac$. In this case, $M_2^n(Lacd)=\{a,d\}$ and $(d,a) \cap Lc=\emptyset$. So the edge that separates $Lcd$ and $Lac$ is deleted in step 2. The rest edges stay unchanged which implies that $M_2^n(Labc) \notin \{\{a,b\},\{b,c\},\{c,a\}\}$. Then we consider the local area formed by $Lab,Lbc,Lcd,Lad,Lbd$ after a flip. From the assumption $M_2^n(Lacd)=\{a,d\}$ and the cyclic order of $a,b,c,d$, we obtain $M_2^n(Labd)=\{a,d\}$. Thus the edge that separates $Lbd$ and $Lab$ is deleted under $\partial_n$. The converse direction, that is considering the inverse flip $W' \longrightarrow W$ is similar. 
    \begin{figure}[H]
    \centering
    \includegraphics[width=13cm]{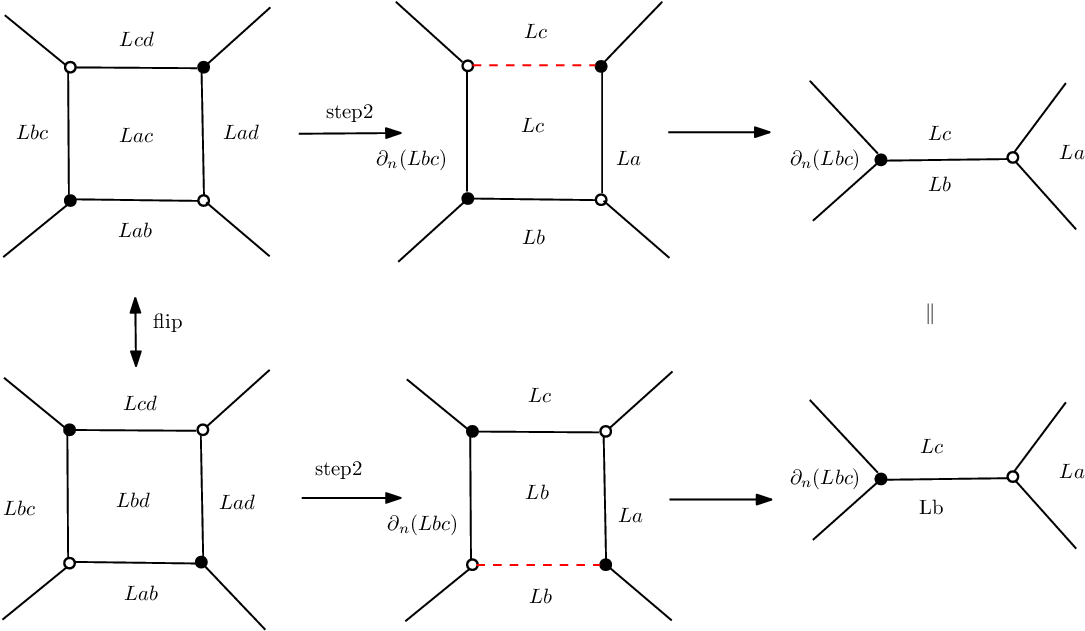}
    \caption{Case 1}
    \label{fig6}
\end{figure}
\noindent Case 2: As in figure \ref{fig7}, two edges are deleted in step 2 in the local area formed by $Lab,Lbc,Lcd,Lad,Lac$. This is equivalent to $M_2^n(Lacd)=\{a,d\}$ with $(d,a) \cap Lc=\emptyset$ and $M_2^n(Labc)=\{a,b\}$ with $(a,b) \cap Lc=\emptyset$. By Lemma \ref{lem1}, the upper left edge outside of the square region must connect to a white point or boundary vertex. So if we consider the boundary map after flip, that is on the local area formed by  $Lab,Lbc,Lcd,Lad,Lbd$ , then edges that separate faces $Lbc,Lcd$ and faces $Lac,Lab$ will be deleted. Conversely, it is also true as case 1.
\begin{figure}[H]
    \centering
    \includegraphics[width=11cm]{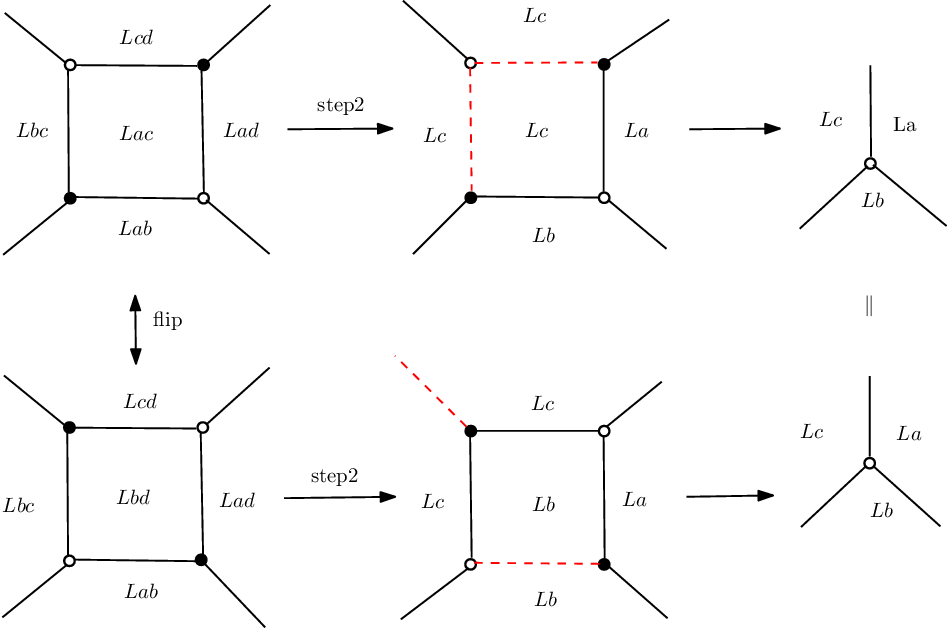}
    \caption{ Case 2}
    \label{fig7}
\end{figure}
Case 3: No edge is deleted in step 2. So the flip $W \longrightarrow W'$ induces a flip $\partial_n(W) \longrightarrow \partial_n(W')$ as in figure \ref{fig9}.
\begin{figure}[H]
    \centering
    \includegraphics[width=11cm]{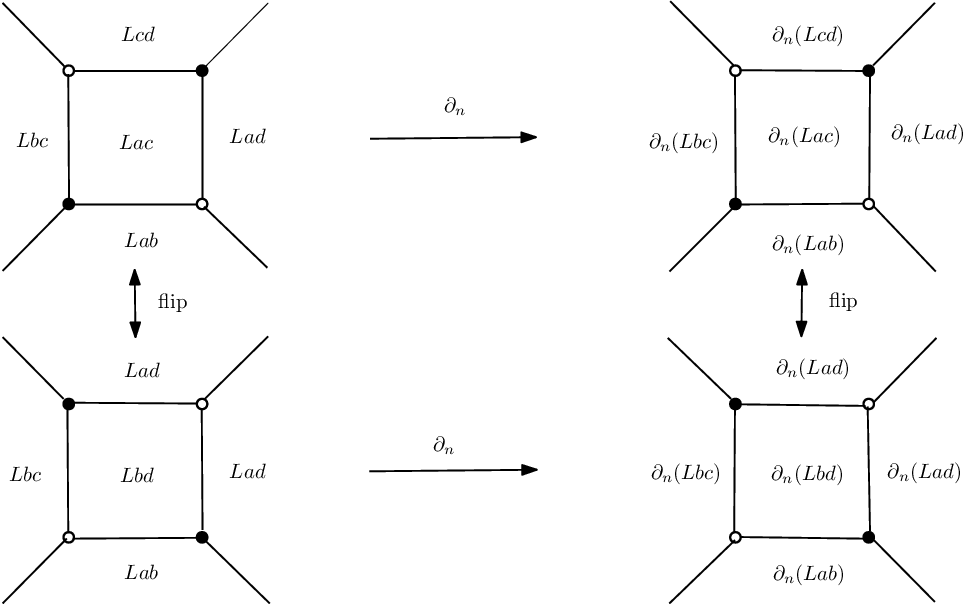}
    \caption{Case 3}
    \label{fig9}
\end{figure}
We have thus shown in cases 1 or 2 that when $\partial_n(W)=\partial_n(W')$, we have $\partial_n(\Sigma_0(W)) \cong \partial_n(\Sigma_0(W'))$. So by the assumption $\partial_n(\Sigma_0(W)) \cong \Sigma_0(\partial_n(W))$, we get $\partial_n(\Sigma_0(W')) \cong \Sigma_0(\partial_n(W'))$. This isomorphism is obvious in case 3. Besides, when $|\partial_n^{-1}|_{W}(\partial_n(J))| \geq 2$, we notice that $\partial_n^{-1}|_{W}(\partial_n(J))$ form a series of consecutive vertices of a white polygon labeled by $\partial_n(J)$ in the plabic tiling of $W$. So the induced subgraph is connected.
\end{proof}

\vspace{0.3cm}
\subsection{Proof of the second theorem and a corollary}\quad
\label{subsec2}

\begin{proof}[Proof of Theorem \ref{thm2}]
$\partial_n(W_0)$ is a maximal w-collection as $W_0=\mathcal{I}_n^k \cup \mathcal{S}_n^k$ by Lemma \ref{lem6}. Suppose that $W$ is obtained from $W_0$ by a sequence of flips, then by Lemma \ref{lem7}, either $\Sigma_0(\partial_n(W_0)) \cong \Sigma_0(\partial_n(W))$ or $\Sigma_0(\partial_n(W))$ is obtained from $\Sigma_0(\partial_n(W_0))$ by a sequence of flips. Since $\Sigma_0(\partial_n(W_0))$ is reduced, then $\Sigma_0(\partial_n(W))$ is also reduced so $\partial_n(W)$ is a maximal w-collection.
\end{proof}
~~~~~~~~~~~~~~~~~~~~~~~~~~~~~~~~~~~~~~~~~~~~~~~~~~~~~~~~~~~~~~~~~~~~~~~~~~~~~~~~~~~~~~~~~~~~~~~~~~~~~~~~~~~~~~~~~~~~~~~~~~~~~~~~~~~~~~~~~~~~~~~~~~~~~~~~~~~~~~~~~~~~~~~~~~~~~~~~~~~~~~~~

\begin{ex}
Take a maximal w-collection $W=\{127,137,136,156,167,135,145,134,123,234,345,456,567\\,678,178,128\}$ in $\binom{[8]}{3}$, then by Definition \ref{def1}, $\partial_8(\Sigma_0(W))$ is produced as showed in Figure \ref{icrafig1}.
    \begin{figure}[H]
    \centering
    \includegraphics[width=13cm]{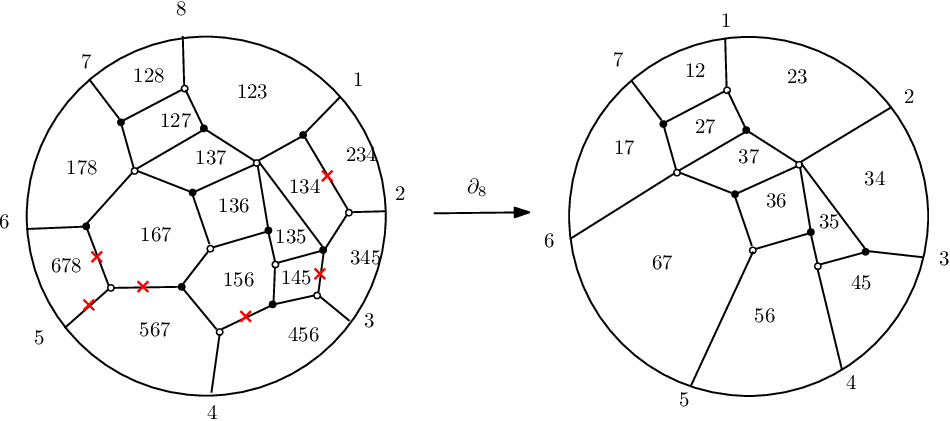}
    \caption{An example of producing reduced plabic graph $\partial_8(\Sigma_0(W))$}
    \label{icrafig1}
\end{figure}
Since $\partial_8(W)=\{27,37,36,35,12,23,34,45,56,67,17\}$ is the collection of face labels on  the reduced plabic graph $\partial_8(\Sigma_0(W))$, it is a maximal w-collection in $\binom{[7]}{2}$
\end{ex}

\vspace{0.3cm}
\begin{cor}
   \label{cor5}
    The translated blades $\{\beta_J | J \in W\}$ of a maximal $w$-collection $W$ induce a finest regular positroid subdivison of $\Delta_{k,n}$.
\end{cor}

\begin{proof}
    Choose any $L \subset [n]$ with $|L|=k-2$, then we use Theorem \ref{thm2} repeatedly. Thus, $\partial_L(W)$ is a maximal w-collection in $\binom{[n] \backslash L}{2}$.  By the same arguments as in the proof of Theorem \ref{thm1}, we obtain that  every octahedral face of $\Delta_{k,n}$ is subdivided. So $\{\beta_J | J \in W\}$ induce a finest positroid subdivision.
\end{proof}

\vspace{0.3cm}

Following Example \ref{exam1},we give an nontrivial example for Corollary \ref{cor5}.
\begin{ex}
     Hypersimplex $\Delta_{3,6}$ is divided into six top dimensional positroid polytopes by blade arrangements $\beta_{135}, \beta_{235},\beta_{145},\beta_{136}$ as showed in Figure \ref{fig17}. For $i \in \{1,2,3,4,5,6\}$, let $\mathcal{M}_i$ be the series-parallel matroid (see \cite{W}) obtained from the spanning tree of a series-parallel graph with labeled edges and let $P_{\mathcal{M}_i}$ be its positroid polytope. Since every series-parallel martoid does not contain an octahedral face (\cite{W}) and the labeled graph does not contain the following forbidden substructure (Figure \ref{fig18}),
    \begin{figure}[H]
    \centering
    \includegraphics[width=8cm]{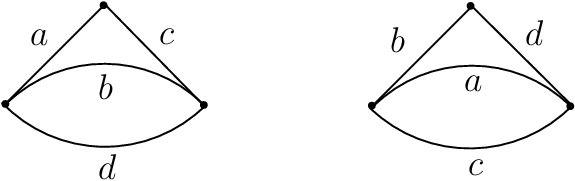}
    \caption{The forbidden structures in the labeled series-parallel graph with $a<b<c<d$ }
    \label{fig18}
\end{figure}

 Then this implies that pyramids $\{Lab,Lbc,Lcd,Lda,Lac\}$ and $\{Lab,Lbc,Lcd,Lda,Lbd\}$ which do not contain non-separated diagonal $\{Lac,Lbd\}$ are not allowed in $\mathcal{M}_i$. Thus this is a finest positroid subdivision.

\begin{figure}[H]
    \centering
    \includegraphics[width=16cm]{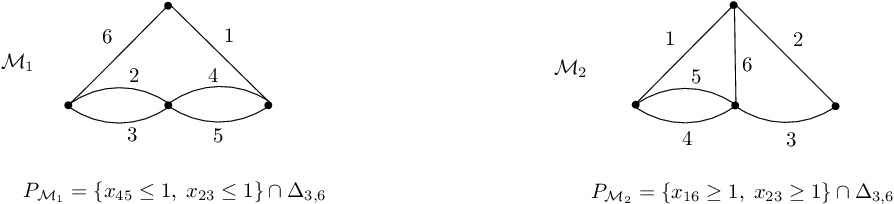}
    \caption*{}
\end{figure}

\begin{figure}[H]
    \centering
    \includegraphics[width=16cm]{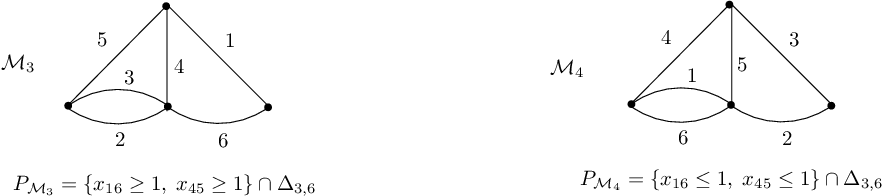}
    \caption*{}
\end{figure}

\begin{figure}[H]
    \centering
    \includegraphics[width=16cm]{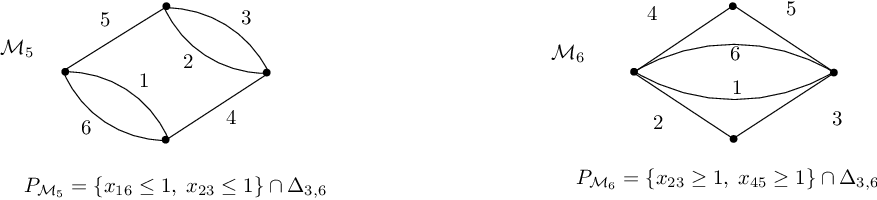}
    \caption{Finest positroid subdivision of $\Delta_{3,6}$ induced by $\beta_{135}, \beta_{235},\beta_{145},\beta_{136}$}
    \label{fig17}
\end{figure}
\end{ex}

\vspace{0.3cm}
\begin{rmk}
    From Corollary \ref{cor5}, we know that there exists unique maximal cone $\mathcal{M}_W$ in $\overline{\mathcal{Z}_{k,n}}$ such that  $\mathcal{P}_W \subseteq \mathcal{M}_W$ for every nonfrozen w-collection $W$. And in Theorem \ref{thm1}, we discussed when $\mathcal{P}_W=\mathcal{M}_W$.
    If two maximal w-collection $W_1$ and $W_2$ are connected by a flip i.e. $W_2=W_1\backslash\{Lac\}\cup \{Lbd\}$, then $\mathcal{M}_{W_1}$ and $\mathcal{M}_{W_2}$ are two adjacent maximal cones in $\overline{\mathcal{Z}_{k,n}}$ and $\mathcal{P}_{W_1} \bigcap \mathcal{P}_{W_2} \subseteq \mathcal{M}_{W_1} \bigcap \mathcal{M}_{W_2}$. We use the Figure \ref{fig14} to present the flip between these objects.
    \begin{figure}[H]
    \centering
    \includegraphics[width=12cm]{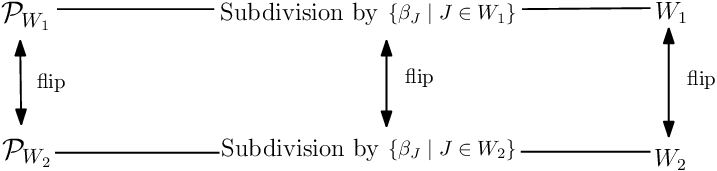}
    \caption{The flip of a positroid subdivision induced by translated blades}
    \label{fig14}
\end{figure}
\end{rmk}

\vspace{1cm}

\section{Conclusions and the further work}

In this paper,we have shown another hypostasis of  maximal weakly separated collection of $\mathcal W_n^k$, they label finest positroids subdivisions of $\Delta_{k,n}$, or the cones in $\overline{\mathcal Z_{k,n}}$.\\
\indent Furthermore, each $W\in \mathcal W_n^k$ is uniquely lifted to a maximal weakly separated collection in $2^{[n]}$ by adding interval sets of size $>k$ and cointerval sets of size $<k$, and is a basis of TP- functions on $2^{[n]}$ satisfying , for any $X$ and $\{i<j<k<l\}\cap X=\emptyset$, 
\[
F(Xj)+F(Xik)=\min ( F(Xi)+F(Xjk), F(Xij)+F(Xk),\]
and 
\[
F(Xik)+F(Xjk)=\min ( F(Xij)+F(Xkl), F(Xil)+F(Xjk).\]

The  supermodular TP-functions, a  subset of which  cut out  by inequalities   $F(Xi)+F(Xj)\le F(Xij)+F(X)$ form a crystal $B(\infty)$  for $SL_n$ \cite{DKK4}.  A supermodular TP-functions is a (sup)  support function to a MV polytope \cite{K}. Let us define cones of supermodular  TP-functions, two supermodular TP-functions $F$ and $G$ belong to the same cone if $F+G$ is a TP-function (it is supermodular since supermodularity is stable under summation). A subdivsion is {\em finest} if all octahedra are subdivided in two halves. Kamnitzer in \cite{K} considers cones in the set of MV-polytopes. MV-polytopes $P$ and $Q$ belong to the same cone if $P+Q$ is an MV-polytope. For example, for $SL_4$ there are 13 such cones of maximal dimension, 12 of which are simplicial with 6 generators and one is not simplicial with 7 generators (Section 6 in \cite{K}).\\
\indent Note that each function defined on vertices on the Boolean cube can be extended as a convex or a concave function to the whole cube $[0,1]^{[n]}$. Thus, we get two dissections of the cube by affinity areas of corresponding extensions. A supermodular TP-function being extended to a convex function on $[0,1]^{[n]}$ is  a (sup) support function  to an MV-polytope. In such a case, the cube is dissected by  Weyl chambers.\\
\indent If we consider a concave extension of a submodular functions, we get a dissection of the cube into generalised polymatroids, since any supermodular TP-function is a $M^{\natural}$-function \cite{M}. One can regard such generalised polymatroids as {\em generalised positroids}.\\
\indent For any $k=2, \cdots, n-2$, the restriction of a supermodular TP-function $F$  to   the vertices of $\Delta_{k,n}$ is a positive tropical Pl\"ucker vector $p$. Hence the intersection of the generalised positroid subdivision  obtained from the concave extension of $F$ with the hyperplane $\sum_ix_i=k$ gives a subdivision  of $\Delta_{k,n}$ into positroids for $p$  (finest subdivison if the subdivison for $F$  is finest).\\
\indent In a subsequent publication we plan to describe cones of finest generalized positroids subdivisions of unit cubes.\\
\indent On the another hand, we can expand a positive tropical Pl\"ucker vector on $\Delta_{k,n}$ to a supermodular TP-function on the Boolean cube $2^{[n]}$ and, hence  get (not uniquely) a generalised positroid subdivision of the cube. 
This extension gives us  finest positroid subdivisions of all $\Delta_{k',n}$,  for all $k'\neq k$, {\em compatible with} $\mathcal D_p$.
\\

\section*{Acknowledgments}
GK thanks School of Mathematical Sciences, Zhejiang University, for hospitality during his visits in 2024. For this project, FL and LZ  are supported by the National Natural Science Foundation of China (No.12131015)


\end{document}